\documentclass[11pt]{amsart}
\usepackage{amsfonts}
\usepackage{amssymb}
\usepackage[letterpaper, margin=1.1in]{geometry}
\usepackage{fancyhdr}

\begin{document}
\author[]{Ali Enayat}
\title[]{The Mostowski Bridge}

\begin{abstract}
In 1950, Novak and Mostowski showed that $\mathsf{GB}$ (G\"{o}del-Bernays
theory of classes) is conservative over $\mathsf{ZF}$, and therefore by G%
\"{o}del's second incompleteness theorem the consistency of $\mathsf{ZF}$ is
unprovable in $\mathsf{GB}$. In the same year Mostowski unveiled a
contrasting result: $\mathsf{GB}$ provides a truth-definition for $\mathsf{ZF%
}$-formulae. Here we first give an expository account of Mostowski's
construction and surrounding results, and then we show that the construction
bridges the domain of Tarski-style truth theories over $\mathsf{PA}$ with
natural extensions of $\mathsf{ACA}_{0}$ (the arithmetical counterpart of $%
\mathsf{GB}$).
\end{abstract}

\maketitle

\begin{center}
\bigskip

\textbf{1. INTRODUCTION}\bigskip
\end{center}

It has long been known that $\mathsf{GB}$ (G\"{o}del--Bernays theory of
classes\footnote{$\mathsf{GB}$ is also referred to as $\mathsf{BG}$, $%
\mathsf{VNB}$ (von Neumann--Bernays) and $\mathsf{NBG}$ (von
Neumann--Bernays--G\"{o}del) in the literature. In some sources $\mathsf{GB}$
includes the global axiom of choice. With all due respect to von Neumann, I
have followed Mostowski's lead in his writings in choosing $\mathsf{GB}$
among these variants.}) is a conservative extension of $\mathsf{ZF}$
(Zermelo--Fraenkel theory of sets), i.e., if $\varphi $ is a sentence in the
language of $\mathsf{ZF}$ that is provable in $\mathsf{GB}$, then $\varphi $
is provable in $\mathsf{ZF}$.\footnote{%
This result is often attributed to Ilse Novak and Andrzej Mostowski. More
specifically, Novak (\cite{Novak (thesis)}, \cite{Novak-FM}) showed that the
consistency of $\mathsf{ZF}$ implies the consistency of $\mathsf{GB}$ by
demonstrating that if $\mathcal{M}\models \mathsf{ZF}$, then $(\mathcal{M},%
\mathfrak{X}_{\mathrm{Def}})\models \mathsf{GB}$, where $\mathfrak{X}_{%
\mathrm{Def}}$ is the collection of $\mathcal{M}$-definable subsets of $M$
(parameters allowed). Mostowski \cite{Mostowski-impredicative} noted that
Novak's method can be used to show (together with the completeness theorem
for first order logic) that $\mathsf{GB}$ is a conservative extension of $%
\mathsf{ZF}$.} This conservativity suggests that $\mathsf{GB}$ is very close
to $\mathsf{ZF}$. In contrast, it is well-known that $\mathsf{ZF}$ is not
finitely axiomatizable, whereas $\mathsf{GB}$ is finitely axiomatizable.
Mostowski \cite{Mostowski-impredicative} came up with an elegant
construction that shows another dramatic difference between $\mathsf{ZF}$
and $\mathsf{GB}$, namely, there is a formula $\mathsf{T}_{\mathsf{Most}}(x)$
-- dubbed the \textit{Mostowski truth predicate} here -- such that for all
sentences $\varphi $ in the language of $\mathsf{ZF}$-set theory, we
have:\medskip

\begin{center}
$\mathsf{GB}\vdash \varphi \leftrightarrow \mathsf{T}_{\mathsf{Most}%
}(\ulcorner \varphi \urcorner ).$\footnote{%
See Corollary 3.6 for a stronger form of this result.}\medskip
\end{center}

\noindent The conservativity of $\mathsf{GB}$ over $\mathsf{ZF}$ combined
with G\"{o}del's second incompleteness theorem makes it clear that $\mathsf{%
Con}(\mathsf{ZF})$\ is unprovable in $\mathsf{GB}$. Together with the above
result about $\mathsf{T}_{\mathsf{Most}}$, this reveals a striking
phenomenon: $\mathsf{GB}$ possesses a truth-predicate for $\mathsf{ZF}$, and
yet the formal consistency of $\mathsf{ZF}$ is unprovable in $\mathsf{GB}$,
which must be due to the lack of sufficient formal induction that is
provable in $\mathsf{GB}$ in order to prove the statement \textquotedblleft
all theorems of $\mathsf{ZF}$ are true\textquotedblright .\footnote{%
Indeed, as shown in Corollary 3.9, $\mathsf{GB}$ cannot even prove the
statement \textquotedblleft all of the \textit{axioms} of $\mathsf{ZF}$ are
deemed true by $\mathsf{T}_{\mathsf{Most}}$\textquotedblright .}\medskip

In this paper we first delve into various aspects of Mostowski's
construction, and then we show that the construction can be used as a bridge
between certain well-known extensions of $\mathsf{ACA}_{0}$ and Tarski-style
truth theories over $\mathsf{PA}$. The paper is organized as follows.
Section 2 contains preliminary definitions and results. Section 3 focuses on
Mostowski's construction and includes many results that came to light since
the appearance of Mostowski's paper \cite{Mostowski-impredicative}. In
Section 4 we use the conceptual framework and machinery developed in Section
3 to link the extensions $\mathsf{ACA}_{0}^{\ast }$ and $\mathsf{ACA}%
_{0}^{\prime }$ of $\mathsf{ACA}_{0}$ with the axiomatic truth theory over $%
\mathsf{PA}$ commonly known as $\mathsf{CT}_{0}[\mathsf{PA}].$ I have
recently explored the analogues of the results in Section 4 in the domain of
set theory, in which certain extensions of $\mathsf{GB}$ are linked with the
set-theoretic counterpart of $\mathsf{CT}_{0}[\mathsf{PA}]$; these results
are presented in \cite{Ali-truth-and-sets}. \medskip

Some of the results of this paper were presented at the Warsaw Logic Seminar
in March 2025. The exposition here has benefited from the collegial comments
of the seminar participants, especially those provided by Leszek Ko\l %
odziejczyk and Mateusz \L e\l yk. I am also grateful to Riki Heck and the
anonymous referee for their constructive suggestions for improving the
paper. The final round of proofreading was carried out with the assistance
of ChatGPT.

\bigskip

\begin{center}
\textbf{2. PRELIMINARIES}\bigskip

\textbf{2.1. Generalities}\medskip
\end{center}

All the theories considered here are first order; some such as $\mathsf{PA}$%
\ and $\mathsf{ZF}$ are one-sorted, and others, such as $\mathsf{GB}$ and $%
\mathsf{ACA}_{0}$, are two-sorted (see below for $\mathsf{GB}$ and $\mathsf{%
ACA}_{0})$. Models of one-sorted theories will be represented using
calligraphic fonts ($\mathcal{M}$, $\mathcal{N}$, etc.) and their universes
will be represented using the corresponding roman fonts ($M$, $N$, etc.).
Models of two-sorted theories will be represented as $(\mathcal{M},\mathfrak{%
X})$, where $\mathcal{M}$ is the structure corresponding to the first sort
(i.e., a model of $\mathsf{PA}$ or $\mathsf{ZF}$), and $\mathfrak{X}$\ is a
collection of subsets of $M$ consisting of objects of the second sort. Given
a theory $U$, we use $\mathcal{L}_{U}$ to refer to the language (signature)
of $U$; for example, $\mathcal{L}_{\mathsf{ZF}}=\{=,\in \}.$\medskip

\noindent Note that $(\mathcal{M},\mathfrak{X})\models \mathsf{GB}$ iff the
following two conditions hold:\medskip

\begin{itemize}
\item If $X_{1},\cdot \cdot \cdot ,X_{n}\in \mathfrak{X}$, then $(\mathcal{M}%
,X_{1},\cdot \cdot \cdot ,X_{n})\models \mathsf{ZF}(X_{1},\cdot \cdot \cdot
,X_{n})$.\footnote{$\mathsf{ZF}(X_{1},\cdot \cdot \cdot ,X_{n})$ denotes the
extension of $\mathsf{ZF}$ in which the names for $X_{1},\cdot \cdot \cdot
,X_{n}$ can be used in the replacement scheme.}

\item If $X_{1},\cdot \cdot \cdot ,X_{n}\in \mathfrak{X},$ and $Y$ is
parametrically definable in $(\mathcal{M},X_{1},\cdot \cdot \cdot ,X_{n})$,
then $Y\in \mathfrak{X}.$\medskip
\end{itemize}

\noindent Similarly, $(\mathcal{M},\mathfrak{X})\models \mathsf{ACA}_{0}$
iff the following two conditions hold:

\begin{itemize}
\item If $X_{1},\cdot \cdot \cdot ,X_{n}\in \mathfrak{X}$, then $(\mathcal{M}%
,X_{1},\cdot \cdot \cdot ,X_{n})\models \mathsf{PA}(X_{1},\cdot \cdot \cdot
,X_{n})$.\footnote{$\mathsf{PA}(X_{1},\cdot \cdot \cdot ,X_{n})$ denotes the
extension of $\mathsf{PA}$ in which the names for the $X_{1},\cdot \cdot
\cdot ,X_{n}$ can be used in the induction scheme.}

\item If $X_{1},\cdot \cdot \cdot ,X_{n}\in \mathfrak{X},$ and $Y$ is
parametrically definable in $(\mathcal{M},X_{1},\cdot \cdot \cdot ,X_{n})$,
then $Y\in \mathfrak{X}.$\medskip
\end{itemize}

\begin{center}
\textbf{2.2. Relationship between Arithmetic and Set Theory}\medskip
\end{center}

\noindent Prima facie $\mathsf{ZF}$ and $\mathsf{PA}$ are unrelated; after
all, $\mathsf{PA}$ axiomatizes basic intuitions concerning the familiar
arithmetical structure of natural numbers, whereas $\mathsf{ZF}$ codifies
the laws governing Cantor's mysterious universe of sets, a vastly more
complex structure. However, thanks to a clever coding idea, introduced first
by Ackermann\textbf{\ }\cite{Ackermann}, $\mathsf{ZF}$ and $\mathsf{PA}$
turn out to be intimately connected at a formal level: $\mathsf{PA}$ is 
\textit{bi-interpretable} with the theory $\mathsf{ZF}^{-\infty }=\mathsf{ZF}%
_{\mathsf{fin}}+\mathsf{TC}$, where $\mathsf{ZF}_{\mathsf{fin}}$ is obtained
from the usual axiomatization of $\mathsf{ZF}$ by replacing the axiom of
infinity with its negation, and $\mathsf{TC}$ expresses \textquotedblleft
the transitive closure of every set exists\textquotedblright .\footnote{%
The role of $\mathsf{TC}$ was elucidated in \cite{ESV} by showing that $%
\mathsf{TC}$ cannot be dropped in this bi-interpretability result.} The
details of this bi-interpretability result are worked out by Kaye and Wong
in \cite{Kaye and Wong}. Their work shows that the two theories are indeed 
\textit{definitionally equivalent} (a.k.a. synonymous), i.e., they have a
common definitional extension. The definitional equivalence of $\mathsf{PA}$
and $\mathsf{ZF}^{-\infty }$, in turn, can be readily employed to show the
definitional equivalence of the theories $\mathsf{ACA}_{0}$ and $\mathsf{GB}%
^{-\infty }$. Here $\mathsf{ACA}_{0}$ is the well-known finitely
axiomatizable subsystem of $\mathsf{Z}_{2}$ (Second Order Arithmetic%
\footnote{%
In the old days $\mathsf{Z}_{2}$ was known as \textit{Formal Analysis}, and $%
\mathsf{PA}$ was known as \textit{Formal\ Number Theory }(and denoted by $%
\mathsf{Z}_{1})$\textit{.} The standard reference for $\mathsf{Z}_{2}$ is
Simpson's monograph \cite{Steve-book}.}), and $\mathsf{GB}^{-\infty }=%
\mathsf{GB}_{\mathsf{fin}}+\mathsf{TC}$, where $\mathsf{GB}_{\mathsf{fin}}$
is obtained from the usual axiomatization of $\mathsf{GB}$ by replacing the
axiom of infinity by its negation. $\mathsf{GB}^{-\infty }$ has not been the
subject of recent work, but it was studied in the 1950s and 1960s; for
example Rabin \cite{Michael} proved that $\mathsf{GB}\backslash \{\mathsf{%
Infinity}\}$ has no recursive model.\footnote{%
Rabin's theorem is one of the precursors of Tennenbaum's theorem on the
nonexistence of nonstandard recursive models of \textsf{PA}.} Also, H\'{a}%
jek and Vop\u{e}nka \cite{Petr-Petr} showed that $\mathsf{TC}$\textsf{\ }is
not provable in the theory $\mathsf{GB}_{\mathsf{fin}}$; later, Hauschild 
\cite{Hauschild} gave a direct construction of a model of $\mathsf{ZF}_{%
\mathsf{fin}}+\lnot \mathsf{TC}$ (Hauschild's construction was rediscovered
and fine-tuned in \cite{ESV}). Another prominent source of study of $\mathsf{%
GB}^{-\infty }$ is Vop\u{e}nka's paper \cite{Petr-1}, in which he showed
that both the power set axiom and the axiom of choice (even global choice)
can be proved from the rest of the axioms of $\mathsf{GB}^{-\infty }$.
\medskip

\begin{itemize}
\item In this paper many definitions and results are formulated in terms of $%
\mathsf{GB}^{\pm \infty }$ (and $\mathsf{ZF}^{\pm \infty })$. Except for the
sentences in this bullet item, all sentences involving the symbol $\pm
\infty $ should be understood as the conjunction of two sentences, one in
which the symbol $\pm $ is replaced by $+$ throughout the sentence,
pertaining to $\mathsf{GB}^{+\infty }=\mathsf{GB,}$ or to $\mathsf{ZF}%
^{+\infty }=\mathsf{ZF}$; and one in which the symbol $\pm $ is replaced by $%
-$ throughout the sentence, pertaining to $\mathsf{GB}^{-\infty }$ (which is
definitionally equivalent to $\mathsf{ACA}_{0}),$ or to $\mathsf{ZF}%
^{-\infty }$ (which is definitionally equivalent to $\mathsf{PA}$).\medskip
\end{itemize}

\noindent The following results are well-known. As mentioned in the
introduction, the first one is due to Novak and Mostowski. The second one is
folklore and can be proved by putting the following ingredients together:
(1) Theorem 2.2.1, (2) G\"{o}del's second incompleteness theorem, (3) the
finite axiomatizability of $\mathsf{GB}^{\mathsf{\pm \infty }}$ (and of $%
\mathsf{ACA}_{0}$), and (4) the fact that $\mathsf{ZF}^{\pm \infty }$ (and
also $\mathsf{PA}$) prove the formal consistency of each of their finite
subtheories.\medskip

\noindent \textbf{2.2.1}.~\textbf{Theorem.}~$\mathsf{GB}^{\mathsf{\pm \infty 
}}$ \textit{is conservative over} $\mathsf{ZF}^{\pm \infty }$; \textit{and} $%
\mathsf{ACA}_{0}$ \textit{is conservative over} $\mathsf{PA}$.\medskip

\noindent \textbf{2.2.2}.~\textbf{Theorem.}~$\mathsf{GB}^{\mathsf{\pm \infty 
}}$ \textit{is not interpretable in }$\mathsf{ZF}^{\pm \infty }$; \textit{and%
} $\mathsf{ACA}_{0}$ \textit{is not interpretable in }$\mathsf{PA}$.\medskip

\begin{center}
\textbf{2.3. Truth and Arithmetic}\medskip
\end{center}

\noindent The theory $\mathsf{CT}^{-}\left( \mathsf{F}\right) $ defined
below is formulated in an \textit{expansion} of $\mathcal{L}_{\mathsf{PA}}$
by adding a fresh \textit{unary} predicate $\mathsf{T}(x)$ ($\mathsf{T}$ for
`true') and a fresh unary predicate $\mathsf{F}(x)$ (note that here $\mathsf{%
F}$ designates a collection of formulae, and not `false'). We often write $%
x\in \mathsf{T}$ instead of $\mathsf{T}(x)$, and $x\in \mathsf{F}$ instead
of $\mathsf{F}(x).$\medskip

\noindent \textbf{2.3.1}.~\textbf{Definition.}~In (a) below we use the
following conventions: $\varphi \in \mathsf{FSent}$ expresses
\textquotedblleft $\varphi $ is an $\mathcal{L}_{\mathsf{PA}}$-sentence
obtained by substituting closed terms of $\mathcal{L}_{\mathsf{PA}}$ for
every free variable of a formula in $\mathsf{F}$\textquotedblright ; $\psi
\vartriangleleft \varphi $ expresses \textquotedblleft $\psi $ is an
immediate subformula of $\varphi $\textquotedblright ; $s$ and $t$ range
over closed terms of $\mathcal{L}_{\mathsf{PA}}$; $s^{\circ }$ denotes the
value of the closed term $s$; $\varphi (v)\in {\mathnormal{\mathsf{F}}}%
^{^{1}}$ expresses \textquotedblleft $\varphi \in \mathsf{F}$ and $\varphi $
has at most one free variable $v$\textquotedblright ; and $\psi \lbrack 
\overset{\cdot }{x}/v]$ is (the code of) the formula obtained by
substituting all occurrences of the variable $v$ in $\psi $ with the numeral
representing $x.$\medskip

\noindent \textbf{(a)} $\mathsf{CT}^{-}\left( \mathsf{F}\right) $ is the
conjunction of $(1)$ through $(5)$ below: \medskip

\begin{enumerate}
\item[$(1)$] $\forall \varphi ,\psi \mathsf{\ }\left[ \left( \varphi \in 
\mathsf{T}\rightarrow \varphi \in \mathsf{FSent}\right) \wedge \left( \left(
\psi \vartriangleleft \varphi \wedge \varphi \in \mathsf{F}\right)
\rightarrow \psi \in \mathsf{F}\right) \right] .$

\item[$(2)$] $\forall s,t\mathsf{\ }\left[ (s=t)\in \mathsf{T}%
\leftrightarrow {s}^{\circ }={t}^{\circ }\right] .$

\item[$(3)$] $\forall \varphi ,\psi \in \mathsf{FSent\ }\left[ \left(
\varphi =\lnot \psi \right) \rightarrow \left( \varphi \in \mathsf{T}%
\leftrightarrow \psi \notin \mathsf{T}\right) \right] \mathsf{.}$

\item[$(4)$] $\forall \varphi ,\psi _{1},\psi _{2}\in \mathsf{FSent\ }\left[
\left( \varphi =\psi _{1}\vee \psi _{2}\right) \rightarrow \left( \varphi
\in \mathsf{T}\leftrightarrow \left( \left( \psi _{1}\in \mathsf{T}\right)
\vee \left( \psi _{2}\in \mathsf{T}\right) \right) \right) \right] \mathsf{.}
$

\item[$(5)$] $\forall \varphi \in \mathsf{FSent},\ \forall \psi (v)\in 
\mathsf{F}^{1}\ \left[ \left( \varphi =\exists v\ \psi (v)\right)
\rightarrow \left( \varphi \in \mathsf{T}\leftrightarrow \exists x\ \psi
\lbrack \overset{\cdot }{x}/v]\in \mathsf{T}\right) \right] .$\medskip
\end{enumerate}

\noindent \textbf{(b)} $\mathsf{CT}^{-}$ is the theory whose axioms are
obtained by substituting the predicate $x\in \mathsf{F}$ by the $\mathcal{L}%
_{\mathsf{PA}}$-formula $x\in \mathsf{Form}_{\mathsf{PA}}$ (expressing
\textquotedblleft $x$ is an $\mathcal{L}_{\mathsf{PA}}$-formula%
\textquotedblright ) in the axioms of $\mathsf{CT}^{-}(\mathsf{F})$. Thus
the axioms of $\mathsf{CT}^{-}$ are formulated in the language $\mathcal{L}_{%
\mathsf{PA}}+\mathsf{T}$ (with no mention of $\mathsf{F}$).\medskip

\noindent \textbf{(c)} $\mathsf{CT}^{-}[\mathsf{PA}]:=\mathsf{PA}+\mathsf{CT}%
^{-}.$\medskip

\noindent \textbf{(d) }\textsf{Int-Ind} (internal induction) is the single
sentence in the language $\mathcal{L}_{\mathsf{PA}}+\mathsf{T}$ that asserts
that every instance of the induction scheme (of \textsf{PA}) is true, i.e., $%
\forall \varphi (x)\in \mathsf{Form}_{\mathsf{PA}}^{1}\ \mathsf{Ind}%
_{\varphi }\in \mathsf{T}$, where $\mathsf{Ind}_{\varphi }$ is the induction
axiom for $\varphi .$\medskip

\noindent \textbf{2.3.2}.~\textbf{Definition.}~Let $\mathcal{M}\models 
\mathsf{PA}$, and suppose $F\subseteq \mathsf{Form}^{\mathcal{M}}$, where $F$
is closed under direct subformulae of $\mathcal{M}$. $\mathsf{FSent}^{(%
\mathcal{M},F)}$ consists of $m\in M$ such that $(\mathcal{M},F)$ satisfies
\textquotedblleft $m$ is an $\mathcal{L}_{\mathsf{PA}}$-sentence obtained by
substituting closed terms of $\mathcal{L}_{\mathsf{PA}}$ for the free
variables of a formula in $\mathsf{F}$\textquotedblright .$\medskip $

\noindent \textbf{(a)} A subset $T$ of $M$ is an $F$-\textit{truth class }on%
\textit{\ }$\mathcal{M}$\textit{\ }if $(\mathcal{M},F,T)\models \mathsf{CT}%
^{-}(\mathsf{F})$; here the interpretation of $\mathsf{F}$ is $F$, and the
interpretation of $\mathsf{T}$ is $T.\medskip $

\noindent \textbf{(b)} A subset $T$ of $M$ is a \textit{full} \textit{truth
class }on\textit{\ }$\mathcal{M}$\textit{\ }if $(\mathcal{M},T)\models 
\mathsf{CT}^{-}$; equivalently: if $(\mathcal{M},F,T)\models \mathsf{CT}^{-}(%
\mathsf{F})$ for $F=\mathsf{Form}^{\mathcal{M}}$.\medskip

\noindent A summary of results about $\mathsf{CT}^{-}[\mathsf{PA}]$ can be
found in \cite{Trio on feasible red.}. For the purposes of this paper, the
following two theorems are especially noteworthy. Theorem 2.3.3 is
essentially due to Kotlarski, Krajewski, and Lachlan \cite{Kotlarski et al}.%
\footnote{%
We are using the qualifier \textquotedblleft essentially\textquotedblright\
since the version of $\mathsf{CT}^{-}[\mathsf{PA}]$ that is used in this
paper is not quite the same as the one used in \cite{Kotlarski et al}. The
conservativity of the version used here was first established by Cie\'{s}li%
\'{n}ski \cite{Cieslinski-book}, using an elaboration of the technique
introduced in \cite{Ali+albert-short}.} Theorem 2.3.4 was first proved in 
\cite{Ali+Albert-long} for the relational version of $\mathsf{CT}^{-}[%
\mathsf{PA}]$; see Leigh \cite{Graham cons.} for the functional version of $%
\mathsf{CT}^{-}[\mathsf{PA}]$. Different proofs for the functional version
of $\mathsf{CT}^{-}[\mathsf{PA]}$ were later presented in \cite{Trio on
feasible red.} and \cite{Ali-curious}.$\medskip $

\noindent \textbf{2.3.3}.~\textbf{Theorem.~}$\mathsf{CT}^{-}[\mathsf{PA}]+$ 
\textsf{Int-Ind} \textit{is conservative over }$\mathsf{PA}$.$\medskip $

\noindent \textbf{2.3.4}.~\textbf{Theorem. }$\mathsf{CT}^{-}[\mathsf{PA}]$%
\textit{\ is interpretable in} $\mathsf{PA}$.$\medskip $

\noindent \textbf{2.3.5}.~\textbf{Remark.}~In contrast to Theorem 2.3.4, $%
\mathsf{CT}^{-}[\mathsf{PA}]+$ \textsf{Int-Ind }is not interpretable in $%
\mathsf{PA.}$ This follows from putting Theorem 2.2.2 with the fact that $%
\mathsf{ACA}_{0}$ is interpretable in $\mathsf{CT}^{-}[\mathsf{PA}]+$ 
\textsf{Int-Ind} (see Lemma 4.3).\medskip

\begin{center}
\textbf{2.4. Truth and Set Theory}\medskip
\end{center}

\noindent \textbf{2.4.1}.~\textbf{Definition.}~Reasoning within $\mathsf{ZF}$%
, for each object $a$ in the universe $\mathrm{V}$ of sets, let $c_{a}$ be a
constant symbol denoting $a$ (e.g., $c_{a}$ can be chosen as the ordered
pair $\left\langle a,0\right\rangle $), and let $\mathcal{L}_{\mathsf{ZF}%
}^{+}$ be the definable proper class of (codes of) sentences of the form $%
\varphi \left( c_{a_{1}},\cdot \cdot \cdot ,c_{a_{n}}\right) $, where $%
\varphi (x_{1},\cdot \cdot \cdot ,x_{n})$ is an $\mathcal{L}_{\mathsf{ZF}}$%
-formula.\medskip

\noindent \textbf{(a)} Suppose $\mathcal{M}$ is a model of $\mathsf{ZF}$. $T$
is an $F$\textit{-truth class} (over $\mathcal{M}$) if $(\mathcal{M},F,T)$
satisfies the conjunction of $(1)$ through $(5)$ below:\medskip

\begin{enumerate}
\item[$(1)$] $\forall x\forall y\left[ \left( \mathsf{T}(x)\rightarrow 
\mathsf{FSent}(x)\right) \wedge \left( y\vartriangleleft x\wedge \mathsf{F}%
(x)\rightarrow \mathsf{F}(y)\right) \right] .$

\item[$(2)$] $\forall a\forall b\left[ \left( \ulcorner c_{a}=c_{b}\urcorner
\in \mathsf{T}\leftrightarrow a=b\right) \wedge \left( \ulcorner c_{a}\in
c_{b}\urcorner \in \mathsf{T}\leftrightarrow a\in b\right) \right] .$

\item[$(3)$] $\forall \varphi ,\psi \in \mathsf{FSent\ }\left[ \left(
\varphi =\lnot \psi \right) \rightarrow \left( \varphi \in \mathsf{T}%
\leftrightarrow \psi \notin \mathsf{T}\right) \right] \mathsf{.}$

\item[$(4)$] $\forall \varphi ,\psi _{1},\psi _{2}\in \mathsf{FSent\ }\left[
\left( \varphi =\psi _{1}\vee \psi _{2}\right) \rightarrow \left( \varphi
\in \mathsf{T}\leftrightarrow \left( \left( \psi _{1}\in \mathsf{T}\right)
\vee \left( \psi _{2}\in \mathsf{T}\right) \right) \right) \right] \mathsf{.}
$

\item[$(5)$] $\forall \varphi \in \mathsf{FSent},\ \forall \psi (v)\in 
\mathsf{F}^{1}\mathsf{\ }\left[ \left( \varphi =\exists v\ \psi (v)\right)
\rightarrow \left( \varphi \in \mathsf{T}\leftrightarrow \exists x\ \mathsf{%
\psi (}c_{x}/v\mathsf{)}\in \mathsf{T}\right) \right] .$\medskip
\end{enumerate}

\noindent \textbf{(b)} In the context of set theory, $\mathsf{CT}^{-}$ is
the theory whose axioms are obtained by substituting the predicate $x\in 
\mathsf{F}$ by the $\mathcal{L}_{\mathsf{ZF}}$-formula $x\in \mathsf{Form}_{%
\mathsf{ZF}}$ (expressing \textquotedblleft $x$ is an $\mathcal{L}_{\mathsf{%
ZF}}$-formula\textquotedblright ) in the axioms of $\mathsf{CT}^{-}(\mathsf{F%
})$. Thus in this context, the axioms of $\mathsf{CT}^{-}$ are formulated in
the language $\mathcal{L}_{\mathsf{ZF}}+\mathsf{T}$ (with no mention of $%
\mathsf{F}$).\medskip

\noindent \textbf{(c)} $\mathsf{CT}^{-}[\mathsf{ZF}]:=\mathsf{ZF}+\mathsf{CT}%
^{-}.$\medskip

\noindent \textbf{(d) }\textsf{Int-Repl} (internal replacement) is the
single sentence in the language $\mathcal{L}_{\mathsf{ZF}}+\mathsf{T}$ that
asserts that every instance of the replacement scheme (of $\mathsf{ZF}$) is
true.\medskip

\noindent \textbf{2.4.2}.~\textbf{Definition.}~Let $\mathcal{M}\models 
\mathsf{ZF}$, and suppose $F\subseteq \mathsf{Form}_{\mathsf{ZF}}^{\mathcal{M%
}}$, where $F$ is closed under direct subformulae of $\mathcal{M}$. $\mathsf{%
FSent}^{(\mathcal{M},F)}$ consists of $m\in M$ such that $(\mathcal{M},F)$
satisfies \textquotedblleft $m$ is an $\mathcal{L}_{\mathsf{ZF}}^{+}$%
-sentence obtained by substituting constants from the class $\left\{
c_{a}:a\in \mathrm{V}\right\} $ for the free variables of a formula in $%
\mathsf{F}$\textquotedblright .$\medskip $

\noindent \textbf{(a)} A subset $T$ of $M$ is an $F$-\textit{truth class }on%
\textit{\ }$\mathcal{M}$\textit{\ }if $(\mathcal{M},F,T)\models \mathsf{CT}%
^{-}(\mathsf{F})$; here the interpretation of $\mathsf{F}$ is $F$ and the
interpretation of $\mathsf{T}$ is $T.$ $\medskip $

\noindent \textbf{(b)} A subset $T$ of $M$ is a \textit{full} \textit{truth
class }on\textit{\ }$\mathcal{M}$\textit{\ }if $(\mathcal{M},T)\models 
\mathsf{CT}^{-}$; equivalently, if $(\mathcal{M},F,T)\models \mathsf{CT}^{-}(%
\mathsf{F})$ for $F=\mathsf{Form}_{\mathsf{ZF}}^{\mathcal{M}}$.\medskip

\noindent Krajewski \cite{Krajewski} showed that $\mathsf{CT}^{-}[\mathsf{ZF}%
]$ is conservative over $\mathsf{ZF}$; indeed his proof can be used to show
the following stronger result in which $\Sigma _{\infty }$-$\mathsf{Sep(T)}$
is the separation scheme for formulae in the language of set theory
augmented with the predicate $\mathsf{T}$; see Theorem 20 of Fujimoto's
paper \cite{Fujimoto-APAL} for a proof.\medskip

\noindent \textbf{2.4.3}.~\textbf{Theorem.}~$\mathsf{CT}^{-}[\mathsf{ZF}%
]+\Sigma _{\infty }$-$\mathsf{Sep(T)}$ \textit{is conservative over} $%
\mathsf{ZF}$.\footnote{$\mathsf{CT}^{-}[\mathsf{ZF}]+\Sigma _{\infty }$-$%
\mathsf{Coll(T)}$ is also conservative over $\mathsf{ZF}$, where $\Sigma
_{\infty }$-$\mathsf{Coll(T)}$ is the collection scheme for formulae in the
language of set theory augmented with the predicate $\mathsf{T}.$ The proof
of this fact will be presented in \cite{Ali-truth-and-sets}.}\medskip

\noindent Since $\Sigma _{\infty }$-$\mathsf{Ind(T)}$ (scheme of induction
over natural numbers in the extended language) is provable in $\mathsf{CT}%
^{-}[\mathsf{ZF}]+\Sigma _{\infty }$-$\mathsf{Sep(T),}$ Theorem 2.4.3
implies:\medskip

\noindent \textbf{2.4.4}.~\textbf{Corollary. }$\mathsf{CT}^{-}\mathsf{[ZF]}$
+ \textquotedblleft $\mathsf{T}$ \textit{is closed under first order proofs}%
\textquotedblright \textit{\ is conservative over }$\mathsf{ZF}$.\medskip

\noindent The above result is in sharp contrast with the well-known fact
that $\mathsf{CT}^{-}\mathsf{[PA]}$ + \textquotedblleft $\mathsf{T}$ is
closed under first order proofs\textquotedblright\ is not conservative over $%
\mathsf{PA}$ since it proves $\mathsf{Con(PA)}$\textsf{; }see Theorem
4.1.\medskip

\noindent The following is a special case of a general result established in 
\cite{Ali+Albert-long}.\medskip

\noindent \textbf{2.4.5}.~\textbf{Theorem.}~$\mathsf{CT}^{-}[\mathsf{ZF}]+$ 
\textsf{Int-Repl} \textit{is conservative over} $\mathsf{ZF}$.\medskip

\noindent \textbf{2.4.6}.~\textbf{Remark. }In contrast with Corollary 2.4.4
and Theorem 2.4.5, the theory $\mathsf{CT}^{-}[\mathsf{ZF}]+$ \textsf{%
Int-Repl} + \textquotedblleft $\mathsf{T}$ is closed under first order proofs%
\textit{\textquotedblright\ }readily implies $\mathsf{Con(ZF)}$\textsf{,}
and is therefore not conservative over $\mathsf{ZF}$. \medskip

\noindent \textbf{2.4.7}.~\textbf{Theorem.}~$\mathsf{CT}^{-}[\mathsf{ZF}%
]+\Sigma _{\infty }$-$\mathsf{Sep(T)}$ \textit{is interpretable in} $\mathsf{%
ZF}$, \textit{but }$\mathsf{CT}^{-}[\mathsf{ZF}]+$ \textsf{Int-Repl} \textit{%
is not interpretable in} $\mathsf{ZF}$.\medskip

\noindent \textbf{Proof.} The proof of Theorem 2.4.3 shows that $\mathsf{CT}%
^{-}[\mathsf{ZF}]+\Sigma _{\infty }$-$\mathsf{Sep(T)}$ is locally
interpretable in $\mathsf{ZF.}$ Since $\mathsf{ZF}$ is a reflective theory
(i.e., proves the formal consistency of each of its finite subtheories), the
global interpretability of $\mathsf{CT}^{-}[\mathsf{ZF}]+\Sigma _{\infty }$-$%
\mathsf{Sep(T)}$ in $\mathsf{ZF}$ then follows by Orey's compactness
theorem. The failure of interpretability of $\mathsf{CT}^{-}[\mathsf{ZF}]+$ 
\textsf{Int-Repl} in $\mathsf{ZF}$ follows from Theorem 2.2.2 and the fact
that $\mathsf{GB}$ is interpretable in $\mathsf{CT}^{-}[\mathsf{ZF}]+$ 
\textsf{Int-Repl} (whose proof is similar to the proof of Lemma 4.3). \hfill 
$\square $\pagebreak

\begin{center}
\textbf{3. THE MOSTOWSKI CONSTRUCTION}\bigskip
\end{center}

\noindent In this section we first present the details of the construction
of the Mostowski truth predicate for $\mathsf{ZF}^{\pm \infty }$ within $%
\mathsf{GB}^{\pm \infty }$, and then we will explore its various
ramifications. Mostowski's original presentation \cite%
{Mostowski-impredicative} focuses on $\mathsf{GB}$, but the construction
works equally well for $\mathsf{GB}^{-\infty }$. In what follows $\omega $
is the collection of finite ordinals as viewed within $\mathsf{GB}^{\pm
\infty }$; note that $\omega $ is a proper class in $\mathsf{GB}^{-\infty }.$%
\medskip

\noindent \textbf{3.1.~Definition.}~Reasoning in $\mathsf{GB}^{\pm \infty }$%
, the \textit{Mostowski cut}, denoted $\mathsf{C}_{\mathsf{Most}}$, consists
of $k\in \omega $ such that there is a class $T$ with the property that $T$
is a $\mathsf{Depth}_{k}$-truth class for the structure $(V,\in )$.\footnote{%
In the terminology of models of arithmetic, a \textit{cut} of a nonstandard
model $\mathcal{M}$ of arithmetic is an initial segment of $\mathcal{M}$
that is closed under immediate successors. This terminology can be readily
applied to $\omega $-nonstandard models of set theory.} Here $V$ is the
class of all sets, and $\mathsf{Depth}_{k}$ is the collection of $\mathcal{L}%
_{\mathsf{ZF}}$-formulae $\varphi $ with $\mathsf{depth}(\varphi )\leq k$;
here $\mathsf{depth}(\varphi )$ is the length of the longest path in the
parsing tree of $\varphi $ (also known as the formation/syntactic tree).
Thus $\mathsf{Depth}_{0}=\varnothing $, and $\mathsf{Depth}_{1}$ consists of
atomic formulae.\footnote{%
We can also define another cut $\mathsf{D}_{\mathsf{Most}}$ consisting of
numbers $k$ such that there is a class $T$ such that $T$ is a $\Sigma _{k}$%
-truth class for the structure $(V,\in )$. All the results of this section
can be established, \textit{mutatis mutandis}, for $\mathsf{D}_{\mathsf{Most}%
}$. Clearly $\mathsf{D}_{\mathsf{Most}}$ is an initial segment of $\mathsf{C}%
_{\mathsf{Most}}.$}\medskip

\noindent \textbf{3.2.~Lemma.}~\textit{Provably in} $\mathsf{GB}^{\pm \infty
}$, $\mathsf{C}_{\mathsf{Most}}$\ \textit{is a cut of }$\omega $.\medskip\ 

\noindent \textbf{Proof.} Suppose $(\mathcal{M},\mathfrak{X})\models \mathsf{%
GB}^{\pm \infty }.$ By default $\mathsf{C}_{\mathsf{Most}}^{(\mathcal{M},%
\mathfrak{X})}$\ contains $0$ since $\mathsf{Depth}_{0}=\varnothing .$ We
will show that $\mathsf{C}_{\mathsf{Most}}^{(\mathcal{M},\mathfrak{X})}$\
contains $1$ and is closed under the successor operation; the definition of%
\textit{\ }$\mathsf{C}_{\mathsf{Most}}$ makes it clear that it is closed
under initial segments. To see that $1\in \mathsf{C}_{\mathsf{Most}}$,
consider the formula $\sigma (x)$ below:\medskip

\begin{center}
$\sigma (x):=\exists a\ \exists b\left[ \left( \left( x=\ulcorner
c_{a}=c_{b}\urcorner \right) \wedge (a=b)\right) \vee \left( \left(
x=\ulcorner c_{a}\in c_{b}\urcorner \right) \wedge (a\in b)\right) \right] .$%
\medskip
\end{center}

\noindent Since $\mathfrak{X}$ includes all $\mathcal{M}$-definable sets,
there is some $T_{1}\in \mathfrak{X}$ such that:\medskip

\begin{center}
$(\mathcal{M},T_{1})\models \forall x(x\in T_{1}\leftrightarrow \sigma (x)).$%
\medskip
\end{center}

\noindent Clearly $T_{1}$ is a $\mathsf{Depth}_{1}$-truth class over $%
\mathcal{M}$. To see that $\mathsf{C}_{\mathsf{Most}}$ is closed under
successors, let $k\in \mathsf{C}_{\mathsf{Most}}$ with $k\geq 1.$ Then by
definition there is some $T\in \mathfrak{X}$\ that is a $\mathsf{Depth}_{k}$%
-truth class over $\mathcal{M}$. Choose $T^{\prime }\in \mathfrak{X}$ such
that:\medskip

\begin{center}
$(\mathcal{M},T,T^{\prime })\models \forall x\ \left[ x\in T^{\prime
}\leftrightarrow \lambda (x,T)\right] ,$ \medskip
\end{center}

\noindent where\medskip

\begin{center}
$\lambda (x,T):=\left[ \left( \mathsf{depth}(x)\leq k\right) \wedge x\in T%
\right] \vee \left( \mathsf{depth}(x)=k+1\wedge \lbrack \mathsf{Neg}(x)\vee 
\mathsf{Exist}(x)\vee \mathsf{Disj}(x)]\right) ,$\medskip
\end{center}

\noindent and\medskip

\begin{center}
$\mathsf{Neg}(x):=\exists y\left[ \left( x=\lnot y\right) \wedge \lnot y\in T%
\right] ,$\medskip

$\mathsf{Exist}(x):=\exists y\ \left[ \exists v\left( x=\exists v\
y(v)\right) \wedge \exists a\ y(c_{a})\in T\right] $, and\medskip

$\mathsf{Disj}(x):=\exists y_{1}\exists y_{2}\left[ \left( x=y_{1}\vee
y_{2}\right) \wedge \left( y_{1}\in T\vee y_{2}\in T\right) \right] .$%
\medskip
\end{center}

\noindent Clearly $T^{\prime }$ is a $\mathsf{Depth}_{k+1}$-truth class over 
$\mathcal{M}$; thus $k+1\in \mathsf{C}_{\mathsf{Most}}$, as desired.\hfill $%
\square $\medskip

\noindent \textbf{3.3.~Definition.}~The following definitions take place
within $\mathsf{GB}^{\pm \infty }.$\medskip

\noindent \textbf{(a)} $\mathsf{Depth}_{\mathsf{Most}}$ consists of $%
\mathcal{L}_{\mathsf{ZF}}$-formulae $\varphi $ such that $\mathsf{depth}%
(\varphi )\in \mathsf{C}_{\mathsf{Most}}.$\medskip

\noindent \textbf{(b)} The \textit{Mostowski truth} \textit{predicate}%
\footnote{%
This definition of $\mathsf{T}_{\mathsf{Most}}$ is clearly $\Sigma _{1}^{1}$%
. However, there are models of $\mathsf{GB}^{\pm \infty }$ in which $\mathsf{%
T}_{\mathsf{Most}}$ is not $\Pi _{1}^{1}$-definable. This is because of the
fact that if $\mathcal{M}$ is a recursively saturated model of $\mathsf{ZF}%
^{\pm \infty }$, then $\left( \mathcal{M},\mathfrak{X}_{\mathrm{Def}}\right) 
$ satisfies $\Delta _{1}^{1}$-comprehension, where $\mathfrak{X}_{\mathrm{Def%
}}$ is the collection of parametrically definable subsets of $\mathcal{M}$
(see $(a)\Rightarrow (b)$ of \cite[Theorem 3.1]{Ali-BS-for-set-theory} for $%
\mathcal{M}\models \mathsf{ZF}$). This implies that if $\mathcal{M}$ is a
recursively saturated model of $\mathsf{ZF}^{\pm \infty }+\lnot \mathsf{Con}(%
\mathsf{ZF}^{\pm \infty })$, then $\mathsf{T}_{\mathsf{Most}}$ is not $\Pi
_{1}^{1}$-definable in the model $\left( \mathcal{M},\mathfrak{X}_{\mathrm{%
Def}}\right) $ of $\mathsf{GB}^{\pm \infty }$.}\textit{\ , }denoted $\mathsf{%
T}_{\mathsf{Most}}(x)$, expresses:\medskip

\begin{center}
$x$ is (the code of) an $\mathcal{L}_{\mathsf{ZF}}^{+}$-formula $\varphi
(c_{a_{1}},...,c_{a_{n}})$ and \medskip

$\exists p\geq \mathsf{depth}(\varphi )$ $\exists T$ $[\varphi
(c_{a_{1}},...,c_{a_{n}})\in T$ and $T$ is a $\mathsf{Depth}_{p}$-truth
class over $(\mathrm{V},\in )]$. \medskip
\end{center}

\noindent \textbf{3.4.~Theorem.}~\textit{Provably in} $\mathsf{GB}^{\pm
\infty }$, $\mathsf{T}_{\mathsf{Most}}(x)$ \textit{is an }$\mathsf{F}$-%
\textit{truth class for }$\mathsf{F}=\mathsf{Depth}_{\mathsf{Most}}$.\medskip

\noindent \textbf{Proof.} Suppose $(\mathcal{M},\mathfrak{X})\models \mathsf{%
GB}^{\pm \infty },$ and let \medskip

\begin{center}
$\mathcal{T}=\{T\in \mathfrak{X}$: $T$ is a $\mathsf{Depth}_{k}$-truth class
over $\mathcal{M}$ for some $k\in \mathsf{C}_{\mathsf{Most}}\}$.\medskip
\end{center}

\noindent Let $T=$ $\mathsf{T}_{\mathsf{Most}}^{(\mathcal{M},\mathfrak{X})}.$
By definition of $\mathsf{T}_{\mathsf{Most}}$, \medskip

\noindent (1) $\ \ T=\cup \mathcal{T}.$ \medskip

\noindent It is easy to see that (1) implies that $T$ satisfies the atomic
clause, the disjunction clause, and the existential clause of part (a) of
Definition 2.4.1 for $F=\mathsf{Depth}_{\mathsf{Most}}^{\mathcal{M}}.$ To
verify the negation clause, suppose $\lnot \varphi \in F.$ We need to check
that:\medskip

\noindent (2) $\ \ \lnot \varphi \in T\leftrightarrow \varphi \notin T.$%
\medskip

\noindent The right-to-left direction of (2) is easy to verify using (1),
but the other direction requires some work. For this purpose, suppose to the
contrary that $\lnot \varphi \in T$ and $\varphi \in T$, and let $\mathsf{%
depth}^{\mathcal{M}}(\varphi )=k$. This implies that there is some $\mathsf{%
Depth}_{k+1}$-truth class $T_{1}\in \mathfrak{X}$ and some $\mathsf{Depth}%
_{k}$-truth class $T_{2}\in \mathfrak{X}$ such that:\medskip

\noindent (3) $\ \ \lnot \varphi \in T_{1},$ and \medskip

\noindent (4) $\ \ \varphi \in T_{2}.$\medskip

\noindent The fact that both $T_{1}$ and $T_{2}$ are elements of $\mathfrak{X%
}$ assures us that $(\mathcal{M},T_{1},T_{2})\models \mathsf{ZF}^{\pm \infty
}(T_{1},T_{2})$, which in turn makes it clear via an induction within $(%
\mathcal{M},T_{1},T_{2})$ on the depth of formulae that $T_{1}$ and $T_{2}$
agree on the truth-evaluation of sentences in $\mathsf{Depth}_{k}$.
Therefore by (4), $\varphi \in T_{1}$, which together with (3) contradicts
the fact that $T_{1}$ is a $\mathsf{Depth}_{k+1}$-truth class over $\mathcal{%
M}$. \hfill $\square $\medskip

\noindent \textbf{3.5.~Remark.}~The last step of the above proof shows that
if $(\mathcal{M},\mathfrak{X})\models \mathsf{GB}^{\pm \infty }$ and $%
T_{i}\in \mathfrak{X}$ is a $\mathsf{Depth}_{k_{i}}$-truth class for $i\in
\{1,2\},$ then $\varphi \in T_{1}$ iff $\varphi \in T_{2}$ for all $\varphi
\in \mathsf{Depth}_{\min \{k_{1},k_{2}\}}^{\mathcal{M}}$.\medskip

\noindent \textbf{3.6.~Corollary.}~\textit{If} $(\mathcal{M},\mathfrak{X}%
)\models \mathsf{GB}^{\pm \infty }$, \textit{then for every standard }$%
\mathcal{L}_{\mathsf{ZF}}$-\textit{formula }$\varphi (x_{1},...,x_{n})$ 
\textit{and any sequence} $\left\langle a_{1},...,a_{n}\right\rangle $ 
\textit{of elements of} $\mathcal{M}$, \textit{the following equivalence
holds}:\medskip

\begin{center}
$\mathcal{M}\models \varphi (a_{1},...,a_{n})$ iff $(\mathcal{M},\mathfrak{X}%
)\models \ulcorner \varphi (c_{a_{1}},...,c_{a_{n}})\urcorner \in \mathsf{T}%
_{\mathsf{Most}}.$ \medskip
\end{center}

\noindent \textbf{3.7.~Theorem.}~\textit{Suppose} $(\mathcal{M},\mathfrak{X}%
)\models \mathsf{GB}^{\pm \infty }$ \textit{and} $\mathcal{M}$ \textit{is }$%
\omega $-\textit{nonstandard}. \textit{Each of the following conditions
implies that} $\mathsf{C}_{\mathsf{Most}}^{(\mathcal{M},\mathfrak{X})}$ =%
\textit{\ the standard cut of }$\mathcal{M}$.\medskip

\noindent \textbf{(a)} $\mathfrak{X}=\mathfrak{X}_{\mathsf{Def(}\mathcal{M}%
\mathsf{)}}$, \textit{where} $\mathfrak{X}_{\mathsf{Def(}\mathcal{M}\mathsf{)%
}}$ \textit{is the collection of} $\mathcal{M}$-\textit{definable subsets of}
$M.$\medskip

\noindent \textbf{(b)} $\mathcal{M}$ \textit{is not recursively saturated}.%
\footnote{%
I am indebted to Mateusz \L e\l yk for drawing my attention to this part of
Theorem 3.7.}\medskip

\noindent \textbf{Proof.} To see that part (a) holds, note that by Lemma
3.2, $\mathsf{C}_{\mathsf{Most}}^{(\mathcal{M},\mathfrak{X})}$ includes the
standard cut of $\mathcal{M}$. The proof of (a) is then immediate in light
of the semantic form of Tarski's undefinability of truth theorem, which
rules out the existence of a set $D\in \mathfrak{X}_{\mathsf{Def(}\mathcal{M}%
\mathsf{)}}$ that has the property that for every standard\textit{\ }$%
\mathcal{L}_{\mathsf{ZF}}$-formula\textit{\ }$\varphi (x_{1},...,x_{n})$,
and any sequence $\left\langle a_{1},...,a_{n}\right\rangle $ of elements of 
$\mathcal{M}$, the following equivalence holds:\medskip

\begin{center}
$\mathcal{M}\models \varphi (a_{1},...,a_{n})$ iff $\mathcal{M}\models
\ulcorner \varphi (c_{a_{1}},...,c_{a_{n}})\urcorner \in D.$
\end{center}

\noindent Part (b) follows from the well-known fact that if $\mathcal{M}$ is
an $\omega $-nonstandard model of $\mathsf{ZF}^{\pm \infty }$ and there is
some $F$-truth class $T$ over $\mathcal{M}$ such that $(\mathcal{M},T)$
satisfies the scheme of induction in the extended language, where $F$
contains all standard $\mathcal{L}_{\mathsf{ZF}}$-formulae, then $\mathcal{M}
$\ is recursively saturated (this is established with an overspill argument;
for a similar argument see \cite[Proposition 15.4]{Kaye Text}).\hfill $%
\square $\medskip

\noindent \textbf{3.8.~Remark.}~Let $\mathsf{GB}_{\mathsf{Def}}^{\pm }$ be
the theory of models of $\mathsf{GB}^{\pm \infty }$ that are of the form $(%
\mathcal{M},\mathfrak{X}_{\mathsf{Def(}\mathcal{M}\mathsf{)}})$, where $%
\mathcal{M}\models \mathsf{ZF}^{\pm \infty }$; i.e., $\mathsf{GB}_{\mathsf{%
Def}}^{\pm }$ is the collection of $\mathcal{L}_{\mathsf{GB}}$-sentences $%
\varphi $ such that $\varphi $ holds in all models of the form $(\mathcal{M},%
\mathfrak{X}_{\mathsf{Def(}\mathcal{M}\mathsf{)}})$, where $\mathcal{M}%
\models \mathsf{ZF}^{\pm \infty }.$ Theorem 5.1 of \cite{Ali+Joel} (whose
proof uses part (a) of Theorem 3.7) shows that:\medskip

\begin{center}
$\mathsf{GB}_{\mathsf{Def}}^{\pm }$ is $\Pi _{1}^{1}$-complete. \medskip
\end{center}

\noindent One can similarly use part (b) of Theorem 3.7 to show that:\medskip

\begin{center}
$\mathsf{GB}_{\lnot \mathsf{recsat}}^{\pm }$ is also $\Pi _{1}^{1}$%
-complete, \medskip
\end{center}

\noindent where $\mathsf{GB}_{\lnot \mathsf{recsat}}^{\pm }$ is the theory
of models of $\mathsf{GB}^{\pm \infty }$ that are of the form $(\mathcal{M},%
\mathfrak{X})$, where $\mathcal{M}$ is a model of $\mathsf{ZF}^{\pm \infty }$
that \textit{is not} recursively saturated. \medskip

\noindent Note that, in contrast, if we let $\mathsf{GB}_{\mathsf{recsat}%
}^{\pm }$ be the theory of models of $\mathsf{GB}^{\pm \infty }$ that are of
the form $(\mathcal{M},\mathfrak{X})$, where $\mathcal{M}$ is a model of $%
\mathsf{ZF}^{\pm \infty }$ that \textit{is} recursively saturated, then in
light of the fact that every model of $\mathsf{GB}^{\pm \infty }$ has an
elementary extension to a recursively saturated model, we have: \medskip

\begin{center}
$\mathsf{GB}_{\mathsf{recsat}}^{\pm }=\{\varphi \in \mathcal{L}_{\mathsf{GB}%
}:\mathsf{GB}^{\pm \infty }\vdash \varphi \},$ \medskip
\end{center}

\noindent which makes it clear that $\mathsf{GB}_{\mathsf{recsat}}^{\pm }$
is recursively enumerable.\medskip

\noindent \textbf{3.9.~Corollary.~ }\textit{Assuming the consistency of }$%
\mathsf{GB}^{\pm \infty }$, \textit{the statement }$\forall x(x\in \omega
\rightarrow x\in \mathsf{C}_{\mathsf{Most}})$ \textit{is unprovable in} $%
\mathsf{GB}^{\pm \infty }.$ \textit{Therefore, there is an instance of} $%
\Sigma _{1}^{1}$-$\mathsf{Ind}$ \textit{that is not provable in} $\mathsf{GB}%
^{\pm \infty }.$\medskip

\noindent \textbf{Proof.} $\mathsf{C}_{\mathsf{Most}}$ is $\Sigma _{1}^{1}$%
-definable within $\mathsf{GB}^{\pm \infty }$, and by Theorem 3.7, there are
(many) $\omega $-nonstandard models $(\mathcal{M},\mathfrak{X})\models 
\mathsf{GB}^{\pm \infty }$ such that $\mathsf{C}_{\mathsf{Most}}^{(\mathcal{M%
},\mathfrak{X})}$ is\textit{\ }the standard cut of\textit{\ }$\mathcal{M}$%
.\hfill $\square $\medskip

\noindent \textbf{3.10.~Definition.}~$\theta _{\mathsf{Most}}$ is the $%
\mathcal{L}_{\mathsf{GB}}$-sentence that asserts that every class is
parametrically definable in the sense of $\mathsf{T}_{\mathsf{Most}}$.%
\footnote{%
The formula $\theta _{\mathsf{Most}}$ was not considered in Mostowski's
paper. To my knowledge, this formula was first introduced in \cite%
{Ali+Albert-long} to show that $\mathsf{CT[PA]}$ is bi-interpretable with $%
\mathsf{ACA}+\theta _{\mathsf{Most}}$, where $\mathsf{CT[PA]}$ is the
extension of $\mathsf{CT}^{-}\mathsf{[PA]}$ with the full induction scheme;
and $\mathsf{ACA}$ is the extension of $\mathsf{ACA}_{0}$ with the full
induction scheme. See Remark 4.6 of this paper for a similar result.} More
specifically, $\theta _{\mathsf{Most}}$ expresses: for every class $X$ there
is an $\mathcal{L}_{\mathsf{ZF}}$-formula $\varphi (x,y)$ and some parameter 
$a\in \mathrm{V}$ such that $X=\left\{ b\in \mathrm{V}:\varphi
(c_{b},c_{a})\in \mathsf{T}_{\mathsf{Most}}\right\} .$\medskip

\noindent \textbf{3.11.~Theorem.}~\textit{Suppose }$(\mathcal{M},\mathfrak{X}%
)\models \mathsf{GB}^{\pm \infty }$. \textit{Then we have}:\medskip

\noindent \textbf{(a)} \textit{If} $\mathcal{M}$ \textit{is} $\omega $-%
\textit{standard, and }$(\mathcal{M},\mathfrak{X})\models \theta _{\mathsf{%
Most}}$, \textit{then} $\mathfrak{X}=\mathfrak{X}_{\mathsf{Def(}\mathcal{M}%
\mathsf{)}}.$\medskip

\noindent \textbf{(b)} $(\mathcal{M},\mathfrak{X}_{\mathsf{Def(}\mathcal{M}%
\mathsf{)}})\models \theta _{\mathsf{Most}}.$\medskip

\noindent \textbf{Proof.} Part (a) readily follows from the definition of $%
\mathfrak{X}_{\mathsf{Def(}\mathcal{M}\mathsf{)}}$. To establish part (b),
use part (a) of Theorem 3.7 together with the reasoning of part (a) of this
theorem.\hfill $\square $\medskip\ 

\noindent \textbf{3.12.~Lemma.}~\textit{Provably in }$\mathsf{GB}^{\pm
\infty },$ \textit{every instance of the replacement scheme whose depth is in%
} $\mathsf{C}_{\mathsf{Most}}$ \textit{is deemed true by} $\mathsf{T}_{%
\mathsf{Most}}.$\medskip

\noindent \textbf{Proof. }Suppose $(\mathcal{M},\mathfrak{X})\models \mathsf{%
GB}^{\pm \infty }$. Suppose $\sigma \in \omega ^{\mathcal{M}}$ is a possibly
nonstandard instance of the replacement scheme within $\mathcal{M}$, thus
for some $\varphi (x,y,z)\in \omega ^{\mathcal{M}}$, $\sigma $ says:\medskip

\textquotedblleft If for some parameter $z$, $\varphi (x,y,z)$ is functional
in $x$ (i.e., $\forall x\ \exists !y\ \varphi (x,y,z)$), then for any set $v$%
, there is a set $w$ such that $w=\left\{ y:\exists x\in v\ \varphi
(x,y,z)\right\} $\textquotedblright $,$ \medskip

\noindent Suppose that:\medskip

\noindent (1) $\ \ \mathsf{depth}^{\mathcal{M}}(\sigma )\in \mathsf{C}_{%
\mathsf{Most}}.$ \medskip

\noindent Furthermore, suppose that for some parameter $m\in M$, \medskip

\noindent (2) $\ \ \left[ \forall x\ \exists !y\ \varphi (x,y,c_{m})\right]
\in \mathsf{T}_{\mathsf{Most}}.$ \medskip

\noindent Then by (1) and (2) there is some $T\in \mathfrak{X}$ such
that:\medskip

\noindent (3) $\ \ T$ is a $\mathsf{Depth}_{k}$-truth class over $\mathcal{M}
$ such that $\left[ \forall x\ \exists !y\ \varphi (x,y,c_{m})\right] \in T.$
\medskip

\noindent Since by assumption $(\mathcal{M},\mathfrak{X})\models \mathsf{GB}%
^{\pm \infty },$ by (3), we have:\medskip

\noindent (4) $\ \ $There is a function $F\in \mathfrak{X}$ with $%
F:M\rightarrow M$ such that:\medskip

\begin{center}
$F(a)=b$ iff \ $\varphi (c_{a},c_{b},c_{m})\in T.$\medskip
\end{center}

\noindent The fact that $F\in \mathfrak{X}$ assures us that $(\mathcal{M}%
,F)\models \mathsf{ZF}(F)$, which implies that $(\mathcal{M},F)$
satisfies:\medskip

\begin{center}
\textquotedblleft for any set $v$, there is a set $w$ such that $w=\left\{
y:\exists x\in v\ \left( y=F(x)\right) \right\} $\textquotedblright .
\medskip
\end{center}

\noindent This makes it clear that $\sigma \in T$, and therefore $\sigma \in 
\mathsf{T}_{\mathsf{Most}}^{(\mathcal{M},\mathfrak{X})}.$\hfill $\square $%
\medskip

\noindent \textbf{3.13.~Lemma.}~\textit{The sentence that expresses the
following is provable in }$\mathsf{GB}^{\pm \infty }:$\medskip

\begin{center}
\textit{For all finite sets }$\Phi \subseteq \mathsf{T}_{\mathsf{Most}}$%
\textit{, and all }$\mathcal{L}_{\mathsf{ZF}}^{+}$-\textit{sentences} $\psi $%
,\medskip

\textit{if }$\Phi \vdash \psi $ \textit{and} $\mathsf{depth}(\psi )\in 
\mathsf{C}_{\mathsf{Most}}$, \textit{then} $\psi \in \mathsf{T}_{\mathsf{Most%
}}.$\medskip
\end{center}

\noindent \textbf{Proof. }Suppose $(\mathcal{M},\mathfrak{X})\models \mathsf{%
GB}^{\pm \infty }$, and let $\Phi $ be an $\mathcal{M}$-finite set of $%
\mathcal{L}_{\mathsf{ZF}}^{+}$-sentences (thus $\Phi \in M$, and the
cardinality of $\Phi $ might be a nonstandard element of $\omega ^{\mathcal{M%
}}$); let $\psi $ be a possibly nonstandard $\mathcal{L}_{\mathsf{ZF}}^{+}$%
-sentence in $\mathcal{M}$; and let $\mathsf{Proof}(x,y,z)$ be the $\mathcal{%
L}_{\mathsf{ZF}}$-formula that expresses \textquotedblleft $x$ is a proof of 
$y$ from assumptions in $z$\textquotedblright . Suppose, furthermore,
that:\medskip

\noindent (1) $\ \ (\mathcal{M},\mathfrak{X})\models \Phi \subseteq \mathsf{T%
}_{\mathsf{Most}}\wedge \mathsf{depth}^{\mathcal{M}}(\psi )\in \mathsf{C}_{%
\mathsf{Most}}\wedge \exists \pi \ \mathsf{Proof}(\pi ,\psi ,\Phi )$.\medskip

\noindent By cut-elimination for first order logic applied within $\mathcal{M%
}$, there is some $\pi ^{\ast }$ in $\mathcal{M}$ such that:\medskip

\noindent (2) $\ \ \mathcal{M}\models \mathsf{Proof}(\pi ^{\ast },\psi ,\Phi
)$ $\wedge $ \textquotedblleft $\pi ^{\ast }$ has the subformula
property\textquotedblright .\medskip

\noindent There is some $k\in \mathsf{C}_{\mathsf{Most}}^{(\mathcal{M},%
\mathfrak{X})}$ such that $k$ dominates the $\mathcal{M}$-depth of $\psi $
and the $\mathcal{M}$-depth of each element of $\Phi ,$ since $\mathsf{depth}%
^{\mathcal{M}}(\psi )\in \mathsf{C}_{\mathsf{Most}}$, $\mathsf{depth}^{%
\mathcal{M}}(\varphi )\in \mathsf{C}_{\mathsf{Most}}$ for all $\varphi \in
\Phi ,$ and $\Phi $ is $\mathcal{M}$-finite. Therefore:\medskip

\noindent (3) $\ \ $There is some $T\in \mathfrak{X}$ such that $T$ is a $%
\mathsf{Depth}_{k}$-truth class over $\mathcal{M}$.\medskip

\noindent Since $T\in \mathfrak{X}$ by (3), $(\mathcal{M},T)\models \mathsf{%
ZF}^{\pm \infty }(T)$, and thus $(\mathcal{M},T)\models \Sigma _{\infty }(T)$%
-$\mathsf{Ind}.$ By (2) and the assumption that $\Phi \subseteq \mathsf{T}_{%
\mathsf{Most}}$, this makes it clear that $\psi \in T$, and therefore $(%
\mathcal{M},\mathfrak{X})$ satisfies $\psi \in \mathsf{T}_{\mathsf{Most}}.$%
\hfill $\square $\medskip

\noindent In what follows $\mathsf{Con}(\mathsf{ZF}^{\pm \infty })$ is the
arithmetical statement asserting $\mathsf{ZF}^{\pm \infty }$ is consistent,
where $\mathsf{ZF}$\ is axiomatized as usual using finitely many axioms
together with the replacement scheme.\medskip

\noindent \textbf{3.14.~Theorem.}~(Solovay) \textit{Provably in }$\mathsf{GB}%
^{\pm \infty },$ $\mathsf{Con}(\mathsf{ZF}^{\pm \infty })$ \textit{holds in} 
$\mathsf{C}_{\mathsf{Most}}.$\footnote{%
Lemma 7.3.1 of Pudl\'{a}k's article \cite{Pavel-Handbook} attributes the
fact that $\mathsf{Con(ZF)}$ is `cut interpretable'\ in $\mathsf{GB}$ to
Solovay's unpublished work.}\medskip

\noindent \textbf{Proof. }Suppose $(\mathcal{M},\mathfrak{X})\models \mathsf{%
GB}^{\pm \infty }$, and let $\mathsf{C}_{\mathsf{Most}}$ denote $\mathsf{C}_{%
\mathsf{Most}}$ as calculated within $(\mathcal{M},\mathfrak{X})$. Let $%
\perp $ denote $\exists x(x\neq x)$ and suppose to the contrary that there
is some $\pi \in \mathsf{C}_{\mathsf{Most}}$ such that:\medskip

\noindent (1) $\ \ \mathcal{M}\models \pi $ is the code of a proof of $\perp 
$ from the axioms of $\mathsf{ZF}^{\pm \infty }$.\medskip

\noindent Here we are using the usual coding of proofs by numbers in which
the code of any of the assumptions in a proof is dominated by the code of
the proof itself. Since $\pi \in \mathsf{C}_{\mathsf{Most}}$, there is some $%
k\in \mathsf{C}_{\mathsf{Most}}$ such that for some $\mathcal{M}$-finite set 
$\Phi $ we have: \medskip

\noindent (2) $\ \ \mathcal{M}\models $ \textquotedblleft each sentence in $%
\Phi $ is an axiom of $\mathsf{ZF}^{\pm \infty }$ of depth at most $k$%
\textquotedblright , \medskip

\noindent and\medskip

\noindent (3) $\ \ \mathcal{M}\models$ \textquotedblleft $\pi $ is the code
of a proof of $\perp $ from $\Phi $\textquotedblright $.$\medskip

\noindent Since all of the axioms of the usual axiomatization of $\mathsf{ZF}
$ (in the sense of $\mathcal{M}$) except for the replacement scheme are of
finite depth in the real world, by Corollaries 3.6 and 3.12, $\Phi \subseteq 
\mathsf{T}_{\mathsf{Most}}$, and by assumption, $\pi \in \mathsf{C}_{\mathsf{%
Most}}$, and $\pi $ is the code of a proof of $\perp $ from $\Phi $ (in the
sense of $\mathcal{M}$), by Lemma 3.13, $\perp \,\in \mathsf{T}_{\mathsf{Most%
}}$, which contradicts Theorem 3.4.\hfill $\square $\medskip

\noindent \textbf{3.15.~Theorem. }$\mathsf{GB}^{\pm \infty }+\forall x(x\in
\omega \rightarrow x\in \mathsf{C}_{\mathsf{Most}})\vdash $
\textquotedblleft All theorems of $\mathsf{ZF}^{\pm \infty }$ are in\textbf{%
\ }$\mathsf{T}_{\mathsf{Most}}$\textquotedblright .\medskip

\noindent \textbf{Proof.} This follows from putting Lemmas 3.12 and 3.13
together.\hfill $\square $\medskip

\noindent \textbf{3.16.~Corollary.}~\textit{The following are provable in} $%
\mathsf{GB}^{\pm \infty }+\forall x(x\in \omega \rightarrow x\in \mathsf{C}_{%
\mathsf{Most}})$.\medskip

\noindent \textbf{(a)} \textit{The consistency of} $\mathsf{ZF}^{\pm \infty
}.$\medskip

\noindent \textbf{(b)} \textit{The consistency of} $\mathsf{GB}^{\pm \infty
}.$\medskip

\noindent \textbf{Proof. }(a) is an immediate consequence of Theorem 3.14.
To see that (b) holds, use the fact that $\mathsf{Con}(\mathsf{ZF}^{\pm
\infty })\rightarrow \mathsf{Con}(\mathsf{GB}^{\pm \infty })$ is already
provable in $\mathsf{I\Delta }_{0}+\mathsf{Supexp}$, and in particular it is
provable in $\mathsf{GB}^{\pm \infty }.$\footnote{%
The proof-theoretic demonstration of conservativity of $\mathsf{GB}$ over $%
\mathsf{ZF}$ was first presented by Shoenfield \cite{Shoenfield} using $%
\varepsilon $-calculus machinery (the proof also shows that $\mathsf{GB}%
^{-\infty }$ \ is conservative over $\mathsf{ZF}^{-\infty }$). It is
well-known that one can also use cut-elimination for the same purpose. Such
proof-theoretic demonstrations are formalizable in $\mathsf{I\Delta }_{0}+%
\mathsf{Supexp}$, where $\mathsf{Supexp}$ asserts the totality of the
superexponential function (also known as tetration). In contrast, the
model-theoretic proof of conservativity of $\mathsf{GB}^{\pm \infty }$ over $%
\mathsf{ZF}^{\pm \infty }$ is formalizable in the extension $\mathsf{I}%
\Sigma _{1}$ of $\mathsf{I\Delta }_{0}+\mathsf{Supexp}$; see \cite[Theorem
IV.4.10]{Hajek and Pudlak}.} \hfill $\square $\medskip

\noindent \textbf{3.17.~Remark.} In light of the fact that $\mathsf{Con}(%
\mathsf{ZF}^{\pm \infty })$ is unprovable in $\mathsf{GB}^{\pm \infty }$ (as
noted in the introduction), Corollary 3.16 provides a different proof of
Corollary 3.9.\medskip

\noindent In the corollary below, $\Omega _{1}$ is the arithmetical
statement asserting the totality of the function $x^{\left\lfloor \log
_{2}(x)\right\rfloor }.$\medskip

\noindent \textbf{3.18.~Corollary.}~$\mathsf{I\Delta }_{0}+\Omega _{1}+%
\mathsf{Con(ZF}^{\pm \infty })$ \textit{is interpretable in} $\mathsf{GB}%
^{\pm \infty }.$ \medskip

\noindent \textbf{Proof.} This follows from putting Theorem 3.14 together
with the well-known technique of `shortening cuts' (as in Section 5 of
Chapter V of H\'{a}jek and Pudl\'{a}k's textbook \cite{Hajek and Pudlak})
that allows one to shorten any definable cut $I$ in a model of Robinson's 
\textsf{Q} to a definable cut $J\subseteq I$ such that $\mathsf{I\Delta }%
_{0}+\Omega _{1}$ holds in $J$. \hfill $\square $\medskip

\noindent In contrast with Corollary 3.18, we have:\medskip

\noindent \textbf{3.19.~Theorem.}~$\mathsf{I\Delta }_{0}+\mathsf{Exp}+%
\mathsf{Con(ZF}^{\pm \infty })$ \textit{is not interpretable in} $\mathsf{GB}%
^{\pm \infty }$. \medskip

\noindent \textbf{Proof.} Within $\mathsf{I\Delta }_{0}+\mathsf{Exp}$ one
can define a cut $J$ such that $\mathsf{Supexp}(x)$ is defined by the
formula $J(x):=\exists y(\mathsf{Supexp}(x)=y)$. Thus within $\mathsf{%
I\Delta }_{0}+\mathsf{Exp}+\mathsf{Con(ZF}^{\pm \infty })$, one can prove
that $\mathsf{Con}(\mathsf{GB}^{\pm \infty })$ holds in $J$ by
cut-elimination, and therefore by shortening $J$\ further using the
shortening technique used in the proof of Corollary 3.18, the following
holds, in which $U\vartriangleright V$ means \textquotedblleft $V$ is
interpretable in $U$\textquotedblright .\medskip

\noindent (1) $\ \ \mathsf{I\Delta }_{0}+\mathsf{Exp}+\mathsf{Con(ZF}^{\pm
\infty })\vartriangleright \mathsf{I\Delta }_{0}+\Omega _{1}+\mathsf{Con(GB}%
^{\pm \infty }).$\medskip

\noindent Thus if $\mathsf{GB}^{\pm \infty }$ could interpret $\mathsf{%
I\Delta }_{0}+\mathsf{Exp}+\mathsf{Con(ZF}^{\pm \infty })$, then by (1) we
would have:\medskip

\noindent (2) $\ \ \mathsf{GB}^{\pm \infty }\vartriangleright \mathsf{%
I\Delta }_{0}+\mathsf{Exp}+\mathsf{Con(ZF}^{\pm \infty })\vartriangleright 
\mathsf{I\Delta }_{0}+\Omega _{1}+\mathsf{Con(GB}^{\pm \infty }).$\medskip

\noindent By transitivity of $\vartriangleright ,$ (2) implies that $\mathsf{%
GB}^{\pm \infty }\vartriangleright \mathsf{I\Delta }_{0}+\Omega _{1}+\mathsf{%
Con(GB}^{\pm \infty })$, which contradicts Pudl\'{a}k's interpretability
form of G\"{o}del's second incompleteness theorem, which states that there
is no consistent recursively enumerable theory such that $U\vartriangleright 
\mathsf{Q}+\mathsf{Con}(U)$ (see Visser's \cite{Visser-Coord-free-G2} for an
exposition). \hfill $\square $\bigskip

\begin{center}
\textbf{4. BRIDGING }$\mathsf{ACA}_{0}$ \textbf{AND} $\mathsf{CT}_{0}[%
\mathsf{PA}]$\bigskip
\end{center}

\noindent $\mathsf{CT}_{0}[\mathsf{PA}]$ is the extension of $\mathsf{CT}%
^{-}[\mathsf{PA}]$ with the scheme $\mathsf{\Delta }_{0}$-$\mathsf{Ind}(%
\mathsf{T})$, i.e., $\mathsf{\Delta }_{0}$-instances of induction in the
extended language $\mathcal{L}_{\mathsf{PA}}+\mathsf{T}$. In this section we
use the concepts and machinery developed in the previous section to connect $%
\mathsf{CT}_{0}[\mathsf{PA}]$ to natural extensions of $\mathsf{ACA}_{0}.$
The following remarkable theorem shows the diverse manners in which $\mathsf{%
CT}_{0}[\mathsf{PA}]$ can be axiomatized. The theorem is due to the
collective work of many researchers (Kotlarski, Cie\'{s}li\'{n}ski, \L e\l %
yk, Wcis\l o, Pakhomov, and the author); see \cite{Bartek+Mateusz on CT0}
and \cite{DC paper of trio} for more information. The following
abbreviations are at work in the statement of Theorem 4.1. \medskip

\begin{itemize}
\item \textquotedblleft $\mathsf{T}$\ is\textit{\ }piecewise
coded\textquotedblright\ is the $\mathcal{L}_{\mathsf{PA}}+\mathsf{T}$%
-sentence that asserts: \medskip
\end{itemize}

\begin{center}
$\forall m\ \exists c\left[ \forall x<m(x\in _{\mathrm{Ack}}c\leftrightarrow
x\in \mathsf{T})\right] .$\footnote{%
Here $x\in _{\mathrm{Ack}}c$ holds iff the $x$-th digit of the binary
expansion of $c$ is 1 ($\in _{\mathrm{Ack}}$ is often referred to as \textit{%
the Ackermann} \textit{epsilon}$).$} \medskip
\end{center}

\begin{itemize}
\item $\mathsf{Prov}_{\varnothing }(\varphi )$ is the $\mathcal{L}_{\mathsf{%
PA}}$-sentence asserting that $\varphi $ is an $\mathcal{L}_{\mathsf{PA}}$%
-sentence that is provable in first order logic alone (with no extra axioms,
hence the symbol $\varnothing $\textsf{).}

\item $\mathsf{Prov}_{\mathsf{PA}}(\varphi )$ is the $\mathcal{L}_{\mathsf{PA%
}}$-sentence asserting that $\varphi $ is an $\mathcal{L}_{\mathsf{PA}}$%
-sentence provable in first order logic from the axioms of $\mathsf{PA}$,
where $\mathsf{PA}$ is formalized by adding the induction scheme to an
appropriate finite collection of axioms.

\item $\mathsf{Prov}_{\mathsf{T}}(\varphi )$ is the $\mathcal{L}_{\mathsf{PA}%
}+\mathsf{T}$-sentence asserting that $\varphi $ is an $\mathcal{L}_{\mathsf{%
PA}}$-sentence that is provable in first order logic from the sentences in $%
\mathsf{T.}$

\item $\mathsf{SentPr}_{\mathsf{T}}(\varphi )$ is the $\mathcal{L}_{\mathsf{%
PA}}+\mathsf{T}$-sentence asserting that $\varphi $ is an $\mathcal{L}_{%
\mathsf{PA}}$-sentence that is provable in sentential (a.k.a. propositional)
logic from the sentences in $\mathsf{T}$.

\item $\mathsf{DC}$ (Disjunctive Correctness) is the $\mathcal{L}_{\mathsf{PA%
}}+\mathsf{T}$-sentence asserting that a finite disjunction is true iff one
of its disjuncts is true; and $\mathsf{DC}_{\mathrm{out}}$ is the `half' of $%
\mathsf{DC}$ that asserts that if a finite disjunction is true, then one of
its disjuncts is true.\footnote{%
The other `half' of $\mathsf{DC}$ is known as $\mathsf{DC}_{\mathrm{in}}.$
As shown in \cite{DC paper of trio}, in contrast to $\mathsf{DC}_{\mathrm{out%
}},$ $\mathsf{CT}^{-}[\mathsf{PA}]+\mathsf{DC}_{\mathrm{in}}$ is
conservative over $\mathsf{PA}$.}\medskip
\end{itemize}

\noindent \textbf{4.1.~The Many Faces Theorem.~}\textit{The following are
equivalent over} $\mathsf{CT}^{-}[\mathsf{I\Delta }_{0}+\mathsf{Exp}]$:$%
\medskip $

\noindent \textbf{(a)} $\mathsf{T}$\ \textit{is piecewise coded}.$\medskip $

\noindent \textbf{(b)} $\mathsf{\Delta }_{0}$-$\mathsf{Ind}(\mathsf{T})$.$%
\medskip $

\noindent \textbf{(c)} $\forall \varphi \left( \mathsf{Prov}_{\varnothing
}(\varphi )\rightarrow \mathsf{T}(\varphi )\right) .\medskip $

\noindent \textbf{(d)} $\forall \varphi \left( \mathsf{Prov}_{\mathsf{T}%
}(\varphi )\rightarrow \mathsf{T}(\varphi )\right) .\medskip $

\noindent \textbf{(e)} $\forall \varphi \left( \mathsf{SentPr}_{\mathsf{T}%
}(\varphi )\rightarrow \mathsf{T}(\varphi )\right) .\medskip $

\noindent \textbf{(f)} $\forall \varphi \left( \mathsf{Prov}_{\mathsf{PA}%
}(\varphi )\rightarrow \mathsf{T}(\varphi )\right) .\medskip $

\noindent \textbf{(g)} $\mathsf{DC}$.$\medskip $

\noindent \textbf{(h)} $\mathsf{DC}_{\mathrm{out}}$.$\medskip $

\noindent \textbf{4.2.~Definition. }Suppose $(\mathcal{M},T)\models \mathsf{%
CT}^{-}[\mathsf{PA}]$. \ $\medskip $

\noindent \textbf{(a)} For each unary $\varphi (x)\in \mathsf{Form}_{\mathsf{%
PA}}^{\mathcal{M}}$, let $\varphi ^{T}:=\left\{ m\in M:(\mathcal{M}%
,T)\models \varphi (\overset{\cdot }{m})\in T\right\} .\medskip $

\noindent \textbf{(b)} $\mathfrak{X}_{\mathrm{Def}_{T}\mathrm{(}\mathcal{M}%
\mathrm{)}}$ is the collection of subsets of $M$ that are of the form $%
\varphi ^{T}$ for some unary $\varphi (x)\in \mathsf{Form}_{\mathsf{PA}}^{%
\mathcal{M}}.\medskip $

\noindent \textbf{4.3.~Lemma.}~\textit{Suppose }$(\mathcal{M},T)\models 
\mathsf{CT}^{-}[\mathsf{PA}]$. \textit{The following are equivalent}:$%
\medskip $

\noindent \textbf{(a)} $(\mathcal{M},T)\models \mathsf{Int}$\textsf{-}$%
\mathsf{Ind}.\medskip $

\noindent \textbf{(b)} $(\mathcal{M},\mathfrak{X}_{\mathrm{Def}_{T}\mathrm{(}%
\mathcal{M}\mathrm{)}})\models \mathsf{ACA}_{0}.$ $\medskip $

\noindent \textbf{Proof. }We only sketch the proof. Recall from part (d) of
Definition 2.3.1 that $\mathsf{Int}$\textsf{-}$\mathsf{Ind}$ is an
abbreviation for the sentence $\forall \varphi (x)\in \mathsf{Form}_{\mathsf{%
PA}}^{1}\ \mathsf{Ind}_{\varphi }\in \mathsf{T}$. The direction (b) $%
\Rightarrow $ (a) is straightforward. The other direction follows from the
following two observations:\medskip

\begin{enumerate}
\item[$(1)$] If $(\mathcal{M},T)\models \mathsf{CT}^{-}[\mathsf{PA}]$, then $%
(\mathcal{M},\mathfrak{X}_{\mathrm{Def}_{T}\mathrm{(}\mathcal{M}\mathrm{)}})$
satisfies all of the axioms of $\mathsf{ACA}_{0}$ with the possible
exception of the induction axiom.

\item[$(2)$] For each unary $\varphi (x)\in \mathsf{Form}_{\mathsf{PA}}^{%
\mathcal{M}}$, $(\mathcal{M},T)\models \left( \mathsf{Ind}_{\varphi }\in
T\right) $ iff for $X=$ $\varphi ^{T},(\mathcal{M},X)\models \mathsf{PA}(X).$%
\hfill $\square $\medskip
\end{enumerate}

\begin{itemize}
\item $\mathsf{ACA}_{0}^{\ast }$ stands for $\mathsf{ACA}_{0}+\forall k\
\exists T\ (T$ is a $\mathsf{Depth}_{k}$-truth class).\footnote{%
Note that $\mathsf{ACA}_{0}^{\ast }$ is definitionally equivalent to $%
\mathsf{GB}^{-\infty }+\forall x(x\in \omega \rightarrow x\in \mathsf{C}_{%
\mathsf{Most}}).$} $\mathsf{ACA}_{0}^{\prime }$ stands for $\mathsf{ACA}%
_{0}+\forall X\ \forall k\ \exists T\ (T$ is a $\mathsf{Depth}_{k}(X)$-truth
class), where $\mathsf{Depth}_{k}(X)$ denotes the set of formulae of depth
at most $k$ in the language $\mathcal{L}_{\mathsf{PA}}+X$ (where $X$\ is a
fresh unary predicate). Note that in this context terms do not contribute to
the calculation of depth of formulae, i.e., atomic formulae involving terms
are counted as having depth 1.\footnote{%
An equivalent axiomatization of these two extensions of $\mathsf{ACA}_{0}$
can be given in terms of the recursion-theoretic concept of the jump
operator as follows: $\mathsf{ACA}_{0}^{\ast }$ is axiomatized as $\mathsf{%
ACA}_{0}+\forall k\ \exists S\ (S=0^{(k)})$, and $\mathsf{ACA}_{0}^{^{\prime
}}$ is axiomatized as $\mathsf{ACA}_{0}+\forall k\ \forall X\ \exists S\
(S=X^{(k)})$.}\medskip
\end{itemize}

\noindent \textbf{4.4.~Lemma.}~$\mathsf{ACA}_{0}^{\ast }+\theta _{\mathsf{%
Most}}$\ \textit{and} $\mathsf{ACA}_{0}^{\prime }+\theta _{\mathsf{Most}}$ 
\textit{axiomatize the same theory }(\textit{where\ }$\theta _{\mathsf{Most}%
} $ \textit{is as in Definition} 3.10).\medskip

\noindent \textbf{Proof.} Note that $\mathsf{ACA}_{0}^{\ast }$ is a
subtheory of $\mathsf{ACA}_{0}^{\prime }$, so it suffices to prove that
given $(\mathcal{M},\mathfrak{X})\models \mathsf{ACA}_{0}^{\ast }+\theta _{%
\mathsf{Most}}$, (1) below holds: \medskip

\noindent (1) $\ \ (\mathcal{M},\mathfrak{X})\models \forall X\ \forall k\
\exists T\ (T$ is a $\mathsf{Depth}_{k}(X)$-truth class).\medskip

\noindent Let $T:=\mathsf{T}_{\mathsf{Most}}^{(\mathcal{M},\mathfrak{X})}$.
Note that since $\theta _{\mathsf{Most}}$ holds in $(\mathcal{M},\mathfrak{X}%
)$, $\mathfrak{X}=\mathfrak{X}_{\mathrm{Def}_{T}\mathrm{(}\mathcal{M}\mathrm{%
)}}.$ Together with the assumption of veracity of $\mathsf{ACA}_{0}^{\ast }$
in $(\mathcal{M},\mathfrak{X}),$ we can conclude:\bigskip

\noindent (2) $\ \ (\mathcal{M},\mathfrak{X}_{\mathrm{Def}_{T}\mathrm{(}%
\mathcal{M}\mathrm{)}})\models \forall k\ \exists T\ (T$ is a $\mathsf{Depth}%
_{k}$-truth class).\medskip

\noindent Given $k\in M$ and $X\in $ $\mathfrak{X}_{\mathrm{Def}_{T}\mathrm{(%
}\mathcal{M}\mathrm{)}}$, choose $p\in M$ such that $X=\varphi ^{T}$ for
some $\varphi \in \mathsf{Depth}_{p}^{\mathcal{M}},$ and let $k^{\prime }\in
M$ such that $\mathcal{M}\models k^{\prime }=k+p.$ Note that any formula in $%
\mathsf{Depth}_{k}^{\mathcal{M}}(X)$ can be translated to an equivalent
formula in $\mathsf{Depth}_{k^{\prime }}$ by substituting all occurrences of 
$t\in X$ (where $t$ is a term) by $\varphi (t)$, and renaming variables so
as to avoid unintended clashes. This makes it clear that (2) implies:\medskip

\noindent (3) $\ \ (\mathcal{M},\mathfrak{X}_{\mathrm{Def}_{T}\mathrm{(}%
\mathcal{M}\mathrm{)}})\models \forall X\ \forall k\ \exists T\ (T $ is a $%
\mathsf{Depth}_{k}(X)$-truth class).

\hfill $\square $\medskip

\begin{itemize}
\item In what follows the expression \textquotedblleft $\omega $%
-interpretable\textit{\textquotedblright\ }simply means \textquotedblleft
interpretable using an interpretation that is the identity on the natural
numbers and arithmetical operations\textquotedblright . Thus in an $\omega $%
-interpretation, numbers are interpreted as numbers, and the interpretation
of each arithmetical operation is itself.\textit{\ }We will write $%
U\vartriangleright _{\omega }V$ to express \textquotedblleft $V$ is $\omega $%
-interpretable in\textit{\ }$U$\textquotedblright .
\end{itemize}

\noindent \textbf{4.5.~Theorem.}~$\mathsf{ACA}_{0}^{\ast }$, $\mathsf{ACA}%
_{0}^{\prime },$ \textit{and} $\mathsf{CT}_{0}[\mathsf{PA}]$ \textit{are
pairwise mutually} $\omega $-\textit{interpretable; hence they agree on
their arithmetical consequences}. \medskip

\noindent \textbf{Proof.} Clearly the identity interpretation witnesses $%
\mathsf{ACA}_{0}^{\prime }\vartriangleright _{\omega }\mathsf{ACA}_{0}^{\ast
}$ since $\mathsf{ACA}_{0}^{\prime }\vdash \mathsf{ACA}_{0}^{\ast }$. So it
suffices to show:\medskip

\begin{center}
$\mathsf{ACA}_{0}^{\ast }\vartriangleright _{\omega }\mathsf{CT}_{0}[\mathsf{%
PA}]\vartriangleright _{\omega }\mathsf{ACA}_{0}^{\prime }.$\medskip
\end{center}

\noindent To see that $\mathsf{ACA}_{0}^{\ast }\vartriangleright _{\omega }%
\mathsf{CT}_{0}[\mathsf{PA}],$\textsf{\ }within $\mathsf{ACA}_{0}^{\ast }$
interpret the truth predicate $\mathsf{T}$ of $\mathsf{CT}_{0}[\mathsf{PA}]$
as $\mathsf{T}_{\mathsf{Most}}.$ To verify that this interpretation
satisfies $\mathsf{CT}_{0}[\mathsf{PA}]$ we will show (1) below, where $T_{%
\mathsf{Most}}$ denotes $\mathsf{T}_{\mathsf{Most}}^{(\mathcal{M},\mathfrak{%
X)}}$.\medskip

\noindent (1) $\ \ $If $(\mathcal{M},\mathfrak{X})\models \mathsf{ACA}%
_{0}^{\ast }$, then $(\mathcal{M},T_{\mathsf{Most}})\models \mathsf{CT}_{0}[%
\mathsf{PA}].$ \medskip

\noindent Suppose $(\mathcal{M},\mathfrak{X})\models \mathsf{ACA}_{0}^{\ast
} $. Then $\mathsf{C}_{\mathsf{Most}}=M,$ and therefore by Theorem 3.4 we
have:\medskip

\noindent (2) $\ \ (\mathcal{M},T_{\mathsf{Most}})\models \mathsf{CT}^{-}[%
\mathsf{PA}].$\medskip

\noindent Thus in light of the equivalence of statements (a) and (b) of
Theorem 4.1 it suffices to show that $T_{\mathsf{Most}}$ is piecewise coded
in $\mathcal{M}$. For this purpose, note that if for some $m\in M,$ and some 
$\varphi \in T_{\mathsf{Most}}$ we have $\varphi <m$, then $\varphi \in 
\mathsf{Depth}_{m}^{\mathcal{M}}$ since we are assuming our coding of
formulae by numbers is a standard one in which the depth of a formula is
dominated by its code. In light of Remark 3.5 (which applies equally well to
this context in light of the bi-interpretation between $\mathsf{ACA}_{0}$
and $\mathsf{GB}^{-\infty }$), we have:\medskip

\noindent (3) $\ \ $If $T\in \mathfrak{X}$ and $T$ is a $\mathsf{Depth}_{m}^{%
\mathcal{M}}$-truth class, then: \medskip

\begin{center}
$\left\{ \varphi \in M:\left( \varphi <^{\mathcal{M}}m\right) \wedge \left(
\varphi \in T_{\mathsf{Most}}\right) \right\} =\left\{ \varphi \in M:\left(
\varphi <^{\mathcal{M}}m\right) \wedge \left( \varphi \in T\right) \right\}
. $\medskip
\end{center}

\noindent Since $(\mathcal{M},T)$ satisfies full induction in the extended
language, $\left\{ \varphi \in M:\left( \varphi <^{\mathcal{M}}m\right)
\wedge \left( \varphi \in T\right) \right\} $ is coded in $\mathcal{M}$ for
each $m\in M$, and therefore by (3) $T_{\mathsf{Most}}$ is piecewise coded
in $\mathcal{M}$. \medskip

Next, we verify that $\mathsf{CT}_{0}[\mathsf{PA}]\vartriangleright _{\omega
}\mathsf{ACA}_{0}^{\prime }\mathsf{.}$ We first describe the interpretation
model-theoretically. Given $(\mathcal{M},T)\models \mathsf{CT}_{0}[\mathsf{PA%
}]$, let $\mathfrak{X}_{\mathrm{Def}_{T}\mathrm{(}\mathcal{M}\mathrm{)}}$ be
the collection of parametrically definable subsets of $\mathcal{M}$ from the
point of view of $T$. More explicitly, using the notation of Definition
4.2:\medskip

\begin{center}
$\mathfrak{X}_{\mathrm{Def}_{T}\mathrm{(}\mathcal{M}\mathrm{)}%
}:=\{X\subseteq M:$ $X=\varphi ^{T}$ for some unary $\varphi (x)\in \mathsf{%
Form}_{\mathsf{PA}}^{\mathcal{M}}\}.$ \medskip
\end{center}

\noindent We claim:\medskip

\noindent (4) $\ \ (\mathcal{M},\mathfrak{X}_{\mathrm{Def}_{T}\mathrm{(}%
\mathcal{M}\mathrm{)}})\models \mathsf{ACA}_{0}^{\prime }.$ \medskip

\noindent It is easy to see that \textsf{Int}-\textsf{Ind} is provable in $%
\mathsf{CT}_{0}[\mathsf{PA}]$, so by Lemma 4.3 we have: \medskip

\noindent (5) $\ \ (\mathcal{M},\mathfrak{X}_{\mathrm{Def}_{T}\mathrm{(}%
\mathcal{M}\mathrm{)}})\models \mathsf{ACA}_{0}.$ \medskip

\noindent Also note that $\theta _{\mathsf{Most}}$ holds in $(\mathcal{M},%
\mathfrak{X}_{\mathrm{Def}_{T}\mathrm{(}\mathcal{M}\mathrm{)}})$. Together
with (5) and Lemma 4.4, to establish (4) it suffices to show that $\mathsf{%
ACA}_{0}^{\ast }$:\medskip

\noindent (6) $\ \ (\mathcal{M},\mathfrak{X}_{\mathrm{Def}_{T}\mathrm{(}%
\mathcal{M}\mathrm{)}})\models \mathsf{ACA}_{0}^{\ast }.$\medskip

\noindent The veracity of (6) follows from the fact, established by Wcis\l o
and \L e\l yk \cite[Lemma 3.7]{Bartek+Mateusz on CT0}, that the following
sentence is provable in $\mathsf{CT}_{0}[\mathsf{PA}]$: \medskip

\noindent (7) $\ \ (\mathcal{M},T)\models \forall k$\ $\exists \varphi \in $ 
$\mathsf{Form}_{\mathsf{PA}}^{1}$ \textquotedblleft $\varphi ^{T}$ is a $%
\mathsf{Depth}_{k}$-truth class\textquotedblright .\medskip\ 

\noindent It is now straightforward to describe the $\omega $-interpretation
of $\mathsf{ACA}_{0}^{\prime }$ within $\mathsf{CT}_{0}[\mathsf{PA}].$
Define an equivalence relation $\thickapprox $ on $\mathsf{Form}_{\mathsf{PA}%
}^{1}$ by:\medskip

\begin{center}
$\varphi _{1}\thickapprox \varphi _{2}$ \ iff $\forall x\left[ \varphi _{1}(%
\overset{\cdot }{x})\in T\leftrightarrow \varphi _{2}(\overset{\cdot }{x}%
)\in T\right] .$\medskip
\end{center}

\noindent Then we define a sort for `sets' within $\mathsf{CT}_{0}[\mathsf{PA%
}]$ whose members are $\mathsf{Form}_{\mathsf{PA}}^{1}$ with equality on
`sets' defined by $\thickapprox $; and we can define the membership $%
\widehat{\in }$ relation between a number $x$ and a `set' $\varphi \in $ $%
\mathsf{Form}_{\mathsf{PA}}^{1}$ via:\medskip

\begin{center}
$x\widehat{\in }\varphi $ \ iff $\ \varphi (\overset{\cdot }{x})\in T.$
\end{center}

\hfill $\square $\medskip

\noindent \textbf{4.6.~Remark.} Let $\theta _{\mathsf{Most}}$ be as in
Definition 3.10. The proof strategy of Theorem 4.5 can be used to show that $%
\mathsf{ACA}_{0}^{\ast }+\theta _{\mathsf{Most}}$ is bi-interpretable with $%
\mathsf{CT}_{0}[\mathsf{PA}].$ Using a similar argument one can show that $%
\mathsf{ACA}$ and $\mathsf{ACA}_{0}+\exists X(X=\mathsf{T}_{\mathsf{Most}})$
are mutually $\omega $-interpretable, and therefore have the same
arithmetical consequences.\medskip

\noindent Theorem 4.7 below was first explicitly stated by McAloon \cite{McA}%
. In the same paper McAloon also presented Theorem 4.9. Proof-theoretic
demonstrations of both Theorems 4.7 and 4.9 were subsequently presented by
Afshari and Rathjen \cite{Afshari-Rathjen-ACA0}, who also pointed out that
Theorem 4.7 is implicit in the work of Jockusch \cite{Jockusch-Ramsey}%
.\medskip

\noindent \textbf{4.7.~Theorem}. $\mathsf{ACA}_{0}^{\prime }$ \textit{and} $%
\mathsf{ACA}_{0}+\mathsf{RT}$ \textit{axiomatize the same theory, where} $%
\mathsf{RT}$ \textit{is the single sentence expressing the Infinite Ramsey
Theorem, i.e.,} \medskip

\begin{center}
$\forall m\ \forall n\ \omega \rightarrow (\omega )_{n}^{m}$. \medskip
\end{center}

\noindent \textbf{4.8.~Definition.~}Given a recursively axiomatized theory $%
U $ extending $\mathrm{I}\Delta _{0}+\mathsf{Exp}$, the \textit{uniform
reflection scheme over} $U$, denoted $\mathsf{REF}(U)$, is defined
via:\medskip

\begin{center}
$\mathsf{REF}(U):=\{\forall x(\mathsf{Prov}_{U}(\ulcorner \varphi (\overset{%
\cdot }{x})\urcorner )\rightarrow \varphi (x)):\varphi (x)\in \mathsf{Form}%
_{U}\}.$\medskip
\end{center}

\noindent The sequence of schemes $\mathsf{REF}^{\alpha }(U)$, where $\alpha
\in \omega $ or $\alpha =\omega $, is defined as follows: \medskip

\noindent $\mathsf{REF}^{0}(U)=U;$\medskip

\noindent $\mathsf{REF}^{n+1}(U)=\mathsf{REF}^{n}(U)+\mathsf{REF}(\mathsf{REF%
}^{n}(U));$\medskip

\noindent $\mathsf{REF}^{\omega }(U)=\bigcup\limits_{n\in \omega }\mathsf{REF%
}^{n}(U).$

\noindent \textbf{4.9.~Theorem.}~\textit{The set of }$\mathcal{L}_{\mathsf{PA%
}}$-\textit{consequences of }$\mathsf{ACA}_{0}^{\prime }$ \textit{is
axiomatized by} $\mathsf{REF}^{\omega }(\mathsf{PA})$.\medskip

\noindent In a parallel development, Kotlarski \cite{Kotlarski bounded}
claimed the following theorem; however, `half' of his proof was subsequently
found to suffer from a nontrivial gap.\footnote{%
Theorem 4.10 asserts that the following holds for all $\mathcal{L}_{\mathsf{%
PA}}$-formulae $\varphi $:
\par
\begin{center}
$\mathsf{CT}_{0}[\mathsf{PA}]\vdash \varphi \Leftrightarrow $ $\mathsf{REF}%
^{\omega }(\mathsf{PA})\vdash \varphi .$%
\end{center}
\par
\noindent The gap in Kotlarski's proof concerns the $(\Leftarrow )$
direction of the above equivalence, but not the other direction (however,
one needs a result of Smory\'{n}ski \cite{Smorynski-Reflection} to show that
the version of $\mathsf{REF}^{\omega }(\mathsf{PA})$ used by Kotlarski is
equivalent to the version used here).} Years later the gap was filled by Wcis%
\l o \cite{Bartek+Mateusz on CT0}, and later by \L e\l yk \cite%
{Mateusz-prolongable} (using a different proof). The latter paper also
includes a model-theoretic proof of $\mathsf{CT}_{0}[\mathsf{PA}]\vdash
\varphi \ \Rightarrow $ $\mathsf{REF}^{\omega }(\mathsf{PA})\vdash \varphi .$
A proof-theoretic demonstration of Theorem 4.10 was presented by Beklemishev
and Pakhomov \cite{Beklemishev-Pakhomov}.\medskip

\noindent \textbf{4.10.~Theorem.}~\textit{The set of }$\mathcal{L}_{\mathsf{%
PA}}$-\textit{consequences of }$\mathsf{CT}_{0}[\mathsf{PA}]$ \textit{is
axiomatized by} $\mathsf{REF}^{\omega }(\mathsf{PA})$.\medskip

\noindent Note that by Theorem 4.5, either of Theorems 4.9 and 4.10 can be
inferred from the other.\medskip

\noindent As mentioned in the introduction, the counterparts of the material
in this section in the realm of set theory have been recently worked out by
the author in the manuscript \cite{Ali-truth-and-sets}, including a `Many
Faces Theorem' for the theory $\mathsf{CT}^{-}[\mathsf{ZF}]\ +$
\textquotedblleft every theorem of \textsf{ZF} is true\textquotedblright .

\noindent \textsc{Ali Enayat, Department of Philosophy, Linguistics, and
Theory of Science, University of Gothenburg}, \textsc{\ Sweden; }\texttt{%
email: ali.enayat@gu.se}


\bigskip

\begin{thebibliography}{\L e\l }
\bibitem[A]{Ackermann} W.~Ackermann,\textit{\ Zur Widerspruchsfreiheit der
Zahlentheorie}, \textbf{Mathematische Annalen,} vol.~117 (1940), pp.~162-194.

\bibitem[AR]{Afshari-Rathjen-ACA0} B.~Afshari and M.~Rathjen, \textit{%
Ordinal Analysis and the Infinite Ramsey Theorem}, in \textbf{How the World
Computes }(edited by\textbf{\ }B.~Cooper et al.), CiE 2012, Lecture Notes in
Computer Science, vol 7318, Springer, Berlin, Heidelberg, 2012, pp.~1-10.

\bibitem[BP]{Beklemishev-Pakhomov} L.~D.\textbf{~}Beklemishev and F.~N.%
\textbf{~}Pakhomov, \textit{Reflection algebras and conservation results for
theories of iterated truth,} \textbf{Annals of Pure and Applied Logic,}
vol.~173(5), 103093 (2022), 41 pp.

\bibitem[C]{Cieslinski-book} C.~Cie\'{s}li\'{n}ski, \textbf{The Epistemic
Lightness of Truth}. \textbf{Deflationism and its Logic}, Cambridge
University Press, Cambridge, 2017.

\bibitem[C\L W]{DC paper of trio} C.~Cie\'{s}li\'{n}ski, M.~\L e\l yk, and
B.~Wcis\l o, \textit{The two halves of disjunctive correctness}, \textbf{%
Journal of Mathematical Logic}, vol.~23(2), Article no.~2250026, (2023), 28
pp.

\bibitem[E-1]{Ali-BS-for-set-theory} A.~Enayat, \textit{Set theoretical
analogues of the Barwise-Schlipf Theorem}, \textbf{Annals of Pure and
Applied Logic}, vol.~173, Paper No. 103158, 16 pages, 2022.

\bibitem[E-2]{Ali-curious} A.~Enayat, \textit{Satisfaction classes with
approximate disjunctive correctness}, \textbf{Review of Symbolic Logic},
vol.~18 (2025), pp.~545-562.

\bibitem[E-3]{Ali-truth-and-sets} A.~Enayat, \textit{Tarskian truth theories
over set theory, }arXiv:2604.03825 [math.LO].

\bibitem[EH]{Ali+Joel} A.~Enayat and J.D.~Hamkins, ZFC \textit{proves that }%
Ord\textit{\ is not weakly compact for definable classes}, \textbf{Journal
of Symbolic Logic}, vol.~83 (2018), pp.~146-164.

\bibitem[E\L W]{Trio on feasible red.} A.~Enayat, M. \L e\l yk, and B. Wcis%
\l o, \textit{Truth and feasible reducibility}, \textbf{Journal of Symbolic
Logic}, vol.~85 (2020), pp.~367-421.

\bibitem[ESV]{ESV} A. Enayat, J.~Schmerl, and A.\textbf{~}Visser, $\omega $%
\textit{-models of finite set theories}, in \textbf{Set Theory, Arithmetic,
and Foundations of Mathematics: Theorems, Philosophies} (edited by
J.~Kennedy and R.~Kossak), Cambridge University Press, 2011, pp.~43-65.

\bibitem[EV-1]{Ali+Albert-long} A. Enayat and A. Visser, \textit{Full
satisfaction classes in a general setting}, privately circulated manuscript
(2012).

\bibitem[EV-2]{Ali+albert-short} A. Enayat and A. Visser, \textit{New
constructions of full satisfaction classes}, in \textbf{Unifying the
Philosophy of Truth} (edited by D.~Achourioti et al.), New York: Springer,
pp.~321-325.

\bibitem[HP]{Hajek and Pudlak} P.~H\'{a}jek and P.~Pudl\'{a}k, \textbf{%
Metamathematics of First-Order Arithmetic}, Springer, 1993.

\bibitem[HV]{Petr-Petr} P.~H\'{a}jek and P.~Vop\u{e}nka,\textit{\ \"{U}ber
die G\"{u}ltigkeit des Fundierungsaxioms in speziellen Systemen der
Mengentheorie}, \textbf{Mathematical Logic Quarterly,} vol.~9 (1963),
pp.~235-241.

\bibitem[H]{Hauschild} K.\textbf{~}Hauschild, \textit{Bemerkungen, das
Fundierungsaxiom betreffend}, \textbf{Mathematical Logic Quarterly,} vol.~12
(1966), pp.~51-56.

\bibitem[F]{Fujimoto-APAL} K.\textbf{~}Fujimoto, \textit{Classes and truths
in set theory}, \textbf{Annals of Pure and Applied Logic, }vol.~163, (2012),
pp.\textbf{~}1484-1523.

\bibitem[J]{Jockusch-Ramsey} C. Jockusch, \textit{Ramsey's theorem and
recursion theory}, \textbf{Journal of Symbolic Logic}, vol.~37 (1972), pp.%
\textbf{~}268-280.

\bibitem[Ka]{Kaye Text} R.\textbf{~}Kaye, \textbf{Models of Peano Arithmetic}%
, Oxford Logic Guides, Oxford University Press, Oxford, 1991.

\bibitem[KW]{Kaye and Wong} R. Kaye and T. L. Wong, \textit{On
interpretations of arithmetic and set theory}, \textbf{Notre Dame Journal of
Formal Logic}, vol.~48 (2007), pp.~497-510.

\bibitem[Ko]{Kotlarski bounded} H. Kotlarski, \textit{Bounded induction and
satisfaction classes}, \textbf{Mathematical Logic Quarterly,} vol.~32
(1986), pp.~531-544.

\bibitem[KKL]{Kotlarski et al} H.~Kotlarski, S.~Krajewski, and
A.~H.~Lachlan, \textit{Construction of satisfaction classes for nonstandard
models}, \textbf{Canadian Mathematical Bulletin,} vol.~24 (1981),
pp.~283-293.

\bibitem[Kr]{Krajewski} S.\textbf{~}Krajewski, \textit{Nonstandard
satisfaction classes}, in \textbf{Set Theory and Hierarchy Theory: A
Memorial Tribute to Andrzej Mostowski} (edited by W.\textbf{~}Marek et al.),
Lecture Notes in Mathematics, vol.~537, Springer-Verlag, Berlin, 1976,
pp.~121-144.

\bibitem[Lei]{Graham cons.} G.~Leigh, \textit{Conservativity for theories of
compositional truth via cut elimination}. \textbf{Journal of Symbolic Logic}%
, vol.~80 (2015), pp.~845-865.

\bibitem[\L e\l ]{Mateusz-prolongable} M.~\L e\l yk, \textit{Model theory
and proof theory of the global reflection principle}, \textbf{Journal of
Symbolic Logic,} vol.~88 (2023), pp.~738-779.

\bibitem[Ma]{McA} K.~McAloon, \textit{Paris-Harrington incompleteness and
progressions of theories},\textbf{\ }in\textbf{\ Proceedings of Symposia in
Pure Mathematics, }vol.\textbf{~} 42, pp.\textbf{~}447-460, 1985.

\bibitem[Mo]{Mostowski-impredicative} A.~Mostowski, \textit{Some
impredicative definitions in the axiomatic set-theory,} \textbf{Fundamenta
Mathematicae,} vol.~37 (1950), pp.~111-124.

\bibitem[N-1]{Novak (thesis)} Novak, I.~L. \textit{On the Consistency of
Goedel's axioms for class and set theory relative to a weaker set of axioms}%
, Doctoral Dissertation, Radcliffe College, 1948.

\bibitem[N-2]{Novak-FM} Novak, I.~L. \textit{A construction for models of
consistent systems}, \textbf{Fundamenta Mathematicae,} vol.~37 (1950),
pp.~87-110.

\bibitem[P]{Pavel-Handbook} P.~Pudl\'{a}k, \textit{Lengths of proofs, }in 
\textbf{Handbook of Proof Theory}, North-Holland, Amsterdam, 1998,
pp.~547-637.

\bibitem[Ra]{Michael} M.~Rabin, \textit{On recursively enumerable and
arithmetic models of set theory}, \textbf{Journal of Symbolic Logic},
vol.~23 (1958), pp.~408-416.

\bibitem[Sh]{Shoenfield} J.~Shoenfield, \textit{A Relative Consistency Proof}%
, \textbf{Journal of Symbolic Logic,} vol.~19 (1954), pp.~21-28.

\bibitem[Si]{Steve-book} S.~Simpson, \textbf{Subsystems of Second Order
Arithmetic}, Perspectives in Mathematical Logic, Springer-Verlag, Berlin,
1999.

\bibitem[Sm]{Smorynski-Reflection} C.~Smory\'{n}ski, $\omega $\textit{%
-consistency and reflection,} in \textbf{Colloque International de Logique:
Clermont-Ferrand}, 1975, Colloq.~Internat. CNRS, vol.~249 (1977), CNRS,
Paris, pp.~167--181.

\bibitem[Vi]{Visser-Coord-free-G2} A.~Visser, \textit{Can We Make the Second
Incompleteness Theorem Coordinate Free?},\textbf{\ Journal of Logic and
Computation}, vol.~21 (2011), pp.~543-560.

\bibitem[Vo]{Petr-1} P.~Vop\u{e}nka, \textit{Axiome der Theorie endlicher
Mengen}, \textbf{\v{C}as.~P\v{e}st.~Mat}. vol.~89 (1964), pp.~312-317.

\bibitem[W\L ]{Bartek+Mateusz on CT0} B.~Wcis\l o and M.~\L e\l yk, \textit{%
Notes on bounded induction for the compositional truth predicate}, \textbf{%
The Review of Symbolic Logic}, vol.~10 (2017), pp.~455-480.\bigskip
\end{thebibliography}
\end{document}